\documentclass[11pt]{imsart}
\usepackage{amsthm,amsmath}
\usepackage[colorlinks,citecolor=blue,urlcolor=blue]{hyperref}

\usepackage{amsmath,amsfonts,amssymb,amsthm}

\usepackage{enumitem}

\def\bbR{\mathbb{R}}
\RequirePackage{amsthm,amsmath,amsfonts,amssymb}
\RequirePackage[numbers]{natbib}
\RequirePackage[colorlinks,citecolor=blue,urlcolor=blue]{hyperref}
\RequirePackage{graphicx}
\usepackage{lineno,hyperref}

\startlocaldefs
\modulolinenumbers[5]

\usepackage[utf8]{inputenc}
\usepackage[english]{babel}

\usepackage[margin=1in]{geometry}
\usepackage{bm}
\usepackage{graphicx}
\usepackage{amssymb,amsmath,amsthm,mathrsfs}
\usepackage{epstopdf}
\usepackage{algorithm}
\usepackage{amsfonts,delarray}
\usepackage{algorithm,algcompatible,amsmath}
\usepackage{rotating,color}
\usepackage{subfigure}
\usepackage{multirow}
\usepackage{titlesec,textcase}
\usepackage{longtable}
\usepackage{bbm}
\usepackage{rotating}

\newtheorem{Theorem}{Theorem}[section]
\newtheorem{Lemma}[Theorem]{Lemma}
\newtheorem{Corollary}[Theorem]{Corollary}
\newtheorem{Proposition}[Theorem]{Proposition}

\newtheorem{Definition}[Theorem]{Definition}

\theoremstyle{definition}

\RequirePackage[normalem]{ulem} 

\definecolor{rp}{RGB}{83,54,106}

\def\boxit#1{\vbox{\hrule\hbox{\vrule\kern6pt\vbox{\kern6pt#1\kern6pt}\kern6pt\vrule}\hrule}}

\endlocaldefs

\begin{document}
\begin{frontmatter}
\title{Information-theoretic Limits for Testing Community Structures in Weighted Networks}

\runtitle{Inference for WSBM}
\runauthor{Mingao Yuan and Zuofeng Shang}
\begin{aug}
\author[A]{\fnms{Mingao} \snm{Yuan}\ead[label=e1]{mingao.yuan@ndsu.edu}}
\and
\author[B]{\fnms{Zuofeng} \snm{Shang}\ead[label=e2]{zshang@njit.edu}}
\address[A]{Department of Statistics, 
North Dakota State University;
\printead{e1}}

\address[B]{Department of Mathematical Sciences,
New Jersey Institute of Technology;
\printead{e2}}
\end{aug}

\begin{abstract}
Community detection refers to the problem of clustering the nodes of a network into groups.
Existing inferential methods for community structure mainly focus on unweighted (binary) networks. Many real-world networks are nonetheless weighted and a common practice is to dichotomize a weighted network to an unweighted one which is known to result in information loss. 
Literature on hypothesis testing in the latter situation is still missing.
In this paper, we study the problem of testing the existence of community structure in weighted networks. Our contributions are threefold: (a). We use the (possibly infinite-dimensional) exponential family to model the weights and derive the sharp information-theoretic limit for the existence of consistent test. Within the limit, any test is inconsistent; and beyond the limit, we propose a useful consistent test. (b). Based on the information-theoretic limits, we provide the first formal way to quantify the loss of information incurred by dichotomizing weighted graphs into unweighted graphs in the context of hypothesis testing. (c). We propose several new and practically useful test statistics. Simulation study show that the proposed tests have good performance. Finally, we apply the proposed tests to an animal social network.
\end{abstract}


\end{frontmatter}

\section{Introduction}
\label{S:1}
In recent decades, network data analysis has attracted increasing attention. One of the most important research topics in network data analysis is to infer the underlying network structures. For instance, in ordinary or hypergraphic stochastic block models (SBM), community detection has been extensively studied by \cite{N01, CY06, ABH16, F10, ZLZ11, ZLZ12, ACB13, LR15, BS16, MS16, A14, A17, ALS18}, among others.
Another research direction is to study hypothesis testing problems regarding the existence of community structures which has been recently studied by
\cite{AV14, VA15, BS16, L16, BM17, GL17a, JKL18, YLFS22, YFS21, YS21,YN20}. 
Existing hypothesis testing literature mainly focus on unweighted networks in which an edge between two nodes is either present or absent. 
Real-world networks are nonetheless often weighted in the sense that the observed edges may be weighted by interaction frequency, volume or similarity, etc; see \cite{A14, AJC15, ALS18, XJL20}.
In weighted networks, edges can be either discrete or continuous random variables characterizing the strength of connectivity (see \cite{A14, AJC15,ALS18, XJL20}).
For example, in airport networks, an edge is weighted by the number of airlines between two airports (\cite{CPV07}); in brain-image networks, the weight can represent the strength of association between two brain regions (\cite{RS10, NBB17}).
To the best of our knowledge, literature on hypothesis testing in general weighted networks 
are still missing, with the exception of some empirical or ad hoc studies (\cite{T18,YYS22}).
Moreover, a common practice in analyzing weighted networks is to dichotomize the weighted edges to binary ones based on which existing unweighted network techniques can be applied. According to some experimental findings, e.g., \cite{A14, AJC15, XJL20}, dichotomizing the weighted networks may result in information loss, whose impact on hypothesis testing is still largely unclear.

In this paper, we are interested in the problem of testing community structures 
in general undirected weighted networks in which the distributions of the weighted edges belong to exponential family.
Our contributions can be summarized into threefold.
(a) We use the general (possibly infinite-dimensional) exponential family to model the distribution of network weights and use contiguity theory coupled with second moment method to derive a sharp information-theoretic limit for the existence of a consistent test. Within the limit, any test is proven inconsistent; and beyond the limit, we propose a useful consistent test. (b) We derive a sharp information-theoretic limit for the existence of a consistent test under dichotomized weighted network, which dramatically differs from the one under original weighted network. Such different information-theoretic limits explicitly quantify the impact of information loss incurred by dichotomizing weighted networks in the context of hypothesis testing. 
(c) We propose useful consistent tests that may outperform the existing empirical or ad hoc approaches proposed in \cite{T18,YYS22}; the latter tests were proposed under the unrealistic assumption that the weights follow some single- or two-parameter distributions
whose applications may be restricted.
More specifically, the limits in (a) and (b) can be characterized by different ellipsoids, and the radius of the latter is significantly larger than the former under any dichotomizing scheme.
The proposed consistent tests are based on weighted signed 
long mixture cycles (WSLMC) which are new in literature with satisfactory numerical performances.
Since the exponential family assumed on the weights is more general than 
Bernoulli distributions assumed on classic unweighted networks, derivation of the limits is substantially more challenging. 

This paper is organized as follows. Section \ref{sec:2} provides a sharp information-theoretic limit for the existence of a consistent test (Section \ref{sec:2.1}) and a quantification of information loss when dichotomizing a weighted network (Section \ref{sec:2.2}). 
In Section \ref{sec:2.3}, the results are strengthened to obtain a more explicit description of the inconsistency of the test, compared with random guessing, under an additional differential equation assumption. 
Practical tests are given in Section \ref{practicalsec}. 
Section \ref{sec:4} involves numerical studies based on simulated and real data. 
Technical proofs are deferred to Section \ref{sec:app}.

\subsection{Problem Formulation}
For a positive integer $n$, let $\mathcal{V}=[n]:=\{1,2,\dots,n\}$ denote the set of network nodes. 
Any node $i\in\mathcal{V}$ is assigned, independently and uniformly at random, a label $\sigma_i\in\{\pm\}$. 
Let $P,Q$ be probability distributions over $\bbR$ and $A\in\bbR^{n\times n}$ be a random symmetric matrix with zero diagonal entries. We say that $A$ follows an undirected weighted stochastic block model (WSBM) with weight distributions $P,Q$, 
denoted $A\sim\mathcal{G}(n, P, Q)$, if $A_{ij},1\le i<j\le n$
are independent and satisfying
\begin{equation*}
  A_{ij} \sim
    \begin{cases}
      $P$, & \text{given $\sigma_i=\sigma_j$},\\
      $Q$, & \text{given $\sigma_i\neq\sigma_j$}.\\
    \end{cases}       
\end{equation*}
When $P,Q$ are Bernoulli distributions, WSBM degenerates to the classic unweighted SBM considered by \cite{A17}. 
For convenience, we call $A_{ij}$'s the network weights.
Clearly, $P=Q$ implies that all network weights are equally distributed regardless of whether a pair of node labels are equal or not,
hence, the network is equivalent to one without community structure.
In practice, only $A$ is observable while $P,Q$ are not. We are interested in the problem of testing whether $P=Q$ given $A$, namely, the following hypothesis testing problem:
\begin{equation}\label{hyp0}
    \textrm{$H_0: P=Q$\,\,\,\, versus\,\,\,\, $H_1: P\neq Q$.}
\end{equation}
A statistical test $T(A)$, a real-valued measurable function of $A$, is said to be consistent if it has asymptotic power approaching one as $n$ goes to infinity. Otherwise, the test is said to be inconsistent. One specific aim of this paper is to study under what circumstances there is a consistent test for (\ref{hyp0}) and how to propose an asymptotically powerful test statistic. This problem will be resolved when $P,Q$ belong to exponential family
that covers a broad range of weighted network models.
In the same setting, \cite{A14, AJC15, XJL20} studied community detection problems. Whereas the problem of testing the existence of community structure still remains elusive.

\subsection{Notation}
For a positive integer $k$,
and a vector of non-negative integers $\alpha=(\alpha_1,\dots,\alpha_k)$, define $|\alpha|=\alpha_1+\alpha_2+\dots+\alpha_k$ and $\alpha!=\alpha_1!\cdots\alpha_k!$. For $\boldsymbol x=(x_1,x_2,\dots,x_k)\bbR^k$, denote $\boldsymbol x^{\alpha}=(x_1^{\alpha_1},x_2^{\alpha_2},\dots,x_k^{\alpha_k})$. For a function $f(\boldsymbol x)$, denote $\partial^{\alpha}f(\boldsymbol x)=\frac{\partial^{\alpha}f(x_1,\dots,x_k)}{\partial x_1^{\alpha_1}\partial x_2^{\alpha_2}\dots \partial x_k^{\alpha_k}}$. Let $Df(\boldsymbol x_0)$
and $D^2f(\boldsymbol x_0)$ denote the gradient and Hessian of $f$ evaluated at $\boldsymbol x_0$, respectively. Let $I(E)$ be the indicator function of an event $E$, $\boldsymbol 1=(1,1,\dots,1)^T$ and $\lambda_{max}(B)$ be the largest eigenvalue of a matrix $B$. Let $D_F^k \psi$ denote the $k$th Fr$\Acute{e}$chet derivative of functional $\psi$ and $D_F^k \psi(f)g^k$ denote the value of multi-linear map $D_F^k \psi(f)$ applied to $(g,g,\dots,g)^T$. Let $\|\boldsymbol x\|_2$ denote the Euclidean norm of $\boldsymbol x $.

\section{Sharp Information-theoretic Limits and Quantification of Information Loss}\label{sec:2}

It might be challenging to derive
sharp information-theoretic limits for testing (\ref{hyp0}) without any regularity assumptions on $P,Q$. Throughout, we assume that $P,Q$ belong to exponential family, which has also been adopted by \cite{AJC15} in community detection.
The parameter dimension of the exponential family can be either finite or infinite, in both cases we shall derive sharp information-theoretic limits. Moreover, we provide the first formal way to quantify the loss of information incurred by dichotomizing the weighted network in the context of hypothesis testing.

\subsection{Sharp Information-theoretic Limit}\label{sec:2.1}
\begin{Definition}\label{def}
Let $\Theta\subset\mathbb{R}^m$ be an $m$-dimensional subset and, for $\boldsymbol\theta\in\Theta$, denote its coordinates $\boldsymbol\theta=(\theta_1,\theta_2,\dots,\theta_m)$. The family of
distributions $\{P_{\boldsymbol\theta}, \boldsymbol\theta\in\Theta\}$ is said to be an exponential family if the probability density of $P_{\boldsymbol\theta}$ has a form
\begin{equation}\label{eqn:ef}
f(x;\boldsymbol\theta)=h(x)e^{\sum_{i=1}^m\theta_iT_i(x)-\psi(\boldsymbol\theta)},
\,\,\,\,x\in\bbR,
\end{equation}
where $h(x),  \psi(\boldsymbol\theta), T_i(x), i=1,2,\dots,m$ are known functions.
\end{Definition}
Let $P$ and $Q$ satisfy (\ref{eqn:ef}) with canonical parameters $\boldsymbol\theta_1$ and $\boldsymbol\theta_2$, respectively.  Let $\boldsymbol\tau=(\tau_1,\tau_2,\dots,\tau_m)^T$ and $\boldsymbol d=(d_1,d_2,\dots,d_m)^T$ be vectors of $m$ fixed constants independent of $n$. Alternatively,
we can express $\boldsymbol\theta_1$ and $\boldsymbol\theta_2$ as follows:
\begin{equation}\label{thetav}
\boldsymbol\theta_1=\boldsymbol\tau-\frac{\boldsymbol\tau_{\boldsymbol d}}{\sqrt{n}},\,\,\,\, \boldsymbol\theta_2=\boldsymbol\tau+\frac{\boldsymbol\tau_{\boldsymbol d}}{\sqrt{n}},
\end{equation}
where $\boldsymbol\tau_{\boldsymbol d}=(\tau_1d_1,\tau_2d_2,\dots,\tau_md_m)^T$.
Under (\ref{thetav}), the hypotheses (\ref{hyp0}) can be rewritten as follows:
\begin{equation}\label{exphyp0}
    \textrm{$H_0:\|\boldsymbol d\|_2=0$\,\,\,\, versus\,\,\,\, $H_1:\|\boldsymbol d\|_2\neq0$.}
\end{equation}
Under $H_0$, $\boldsymbol\theta_1=\boldsymbol\theta_2$, and hence, $P=Q$. Under $H_1$, $P$ and $Q$ are different since at least one component of $\boldsymbol d$ is nonzero.

To derive the sharp limit, firstly we propose a novel test statistic for (\ref{exphyp0}) as follows. 
For integer $k=\log\log\log n$ and distinct nodes $i_1,i_2,\dots,i_k$, let $\mathcal{C}(i_1,i_2,\dots,i_k)$ be the set of all the circular permutations of $i_1,i_2,\dots,i_k$. It is well-known that $\mathcal{C}(i_1,i_2,\dots,i_k)$ has $\frac{(k-1)!}{2}$ elements. Let $\mathcal{I}_n=\{(i_1,i_2,\dots,i_k)\in \mathcal{C}(j_1,j_2,\dots,j_k)|1\leq j_1<j_2<\dots<j_k\leq n\}$. Define the weighted-signed-long-mixture-cycle (WSLMC) test statistic $\mathcal{Z}_n$ as follows:
\begin{eqnarray*}
\mathcal{Z}_n=\frac{\sum_{(i_1,i_2,\dots,i_k)\in\mathcal{I}_n}\prod_{t=1}^k\frac{\boldsymbol\tau_{\boldsymbol d}^T\left(\boldsymbol T(A_{i_ti_{t+1}})-D\psi(\boldsymbol\tau)\right)}{\sqrt{n}}}{\sqrt{\frac{1}{2k}\left[\boldsymbol\tau_{\boldsymbol d}^TD^2\psi(\boldsymbol\tau)\boldsymbol\tau_{\boldsymbol d}\right]^k}},
\end{eqnarray*}
where $\boldsymbol T(x)=(T_1(x),T_2(x),\dots,T_m(x))^T$ and $i_{k+1}=i_1$. Note that each circular permutation of $(i_1,i_2,\dots,i_k)\in\mathcal{I}_n$ can be considered as a cycle. In this sense, the numerator of $\mathcal{Z}_n$ just counts the number of some weighted cycles. Here ``signed'' means $\boldsymbol T(x)$ is centered by subtracting its mean $D\psi(\boldsymbol\tau)$; ``weighted'' means $\boldsymbol T(x)-D\psi(\boldsymbol\tau)$ has weight $\boldsymbol\tau_{\boldsymbol d}$; ``long'' means the length $k$ of the circular permutation (or cycle) goes to infinity as $n$ tends to infinity; `mixture' means $\boldsymbol\tau_{\boldsymbol d}^T\left(\boldsymbol T(A_{i_ti_{t+1}})-D\psi(\boldsymbol\tau)\right)$ is a weighted sum of $m$ terms. This test statistic is motivated by but significantly different from the long-cycle test in \cite{MNS15} and signed-cycle test in \cite{BDER16}.  The WSLMC test rejects $H_0$ if $|\mathcal{Z}_n|>C$ for some constant $C$ dependent on the type I error.

The following Theorem \ref{theorem:1} provides the first sharp information-theoretic limit for existence of consistent test in the weighted network case.

\begin{figure}
    \centering
    \includegraphics[height=7cm,width=10cm]{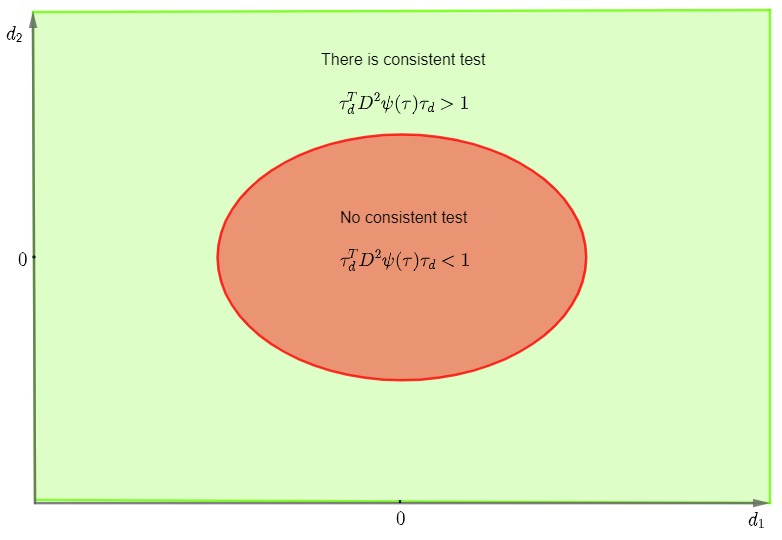}
    \caption{Red region: no consistent test. Green region: there is consistent test.}
    \label{phase}
\end{figure}

\begin{Theorem}\label{theorem:1} Suppose $\partial^{\alpha}\psi$, for $\alpha$ with $|\alpha|=5$, exist and are uniformly bounded for all $\boldsymbol\theta\in\Theta$. Then the following results hold.
\begin{itemize}
\item[(I)] If $\boldsymbol\tau_{\boldsymbol d}D^2\psi(\boldsymbol\tau)\boldsymbol\tau_{\boldsymbol d}<1$, any test is inconsistent. 
\item[(II)] If $\boldsymbol\tau_{\boldsymbol d}D^2\psi(\boldsymbol\tau)\boldsymbol\tau_{\boldsymbol d}>1$, the WSLMC test is consistent.
\end{itemize}
\end{Theorem}
 Theorem \ref{theorem:1} says that, when $\boldsymbol\tau_{\boldsymbol d}^TD^2\psi(\boldsymbol\tau)\boldsymbol\tau_{\boldsymbol d}<1$, any statistical test for (\ref{exphyp0}) cannot achieve asymptotic power one. When $\boldsymbol\tau_{\boldsymbol d}^TD^2\psi(\boldsymbol\tau)\boldsymbol\tau_{\boldsymbol d}>1$, the WSLMC test is consistent and hence optimal in this sense. The solution set  $\{\boldsymbol{d}:\boldsymbol\tau_{\boldsymbol d}^TD^2\psi(\boldsymbol\tau)\boldsymbol\tau_{\boldsymbol d}=1\}$ is an elliptic curve in $\boldsymbol d$, which is the sharp boundary for the existence of a consistent test. 
 When $m=2$, these regions are demonstrated in Figure \ref{phase}.

We point out the WSLMC test is not directly applicable in practice, since the vector $\boldsymbol\tau_{\boldsymbol d}$ and the function $\boldsymbol T(x)=(T_1(x),T_2(x),\dots,T_m(x))^T$ are unknown. Our results only provide some theoretical insights and serve as a benchmark for developing practical statistical tests. Motivated by $\mathcal{Z}_n$, more practical tests shall be proposed in Section \ref{practicalsec}.

Next, we illustrate our results by restricting $P,Q$ to the exponential distribution and the normal distribution.\\

\noindent{\bf Example 1 (Exponential distribution).}
The exponential distribution has density
\begin{equation}\label{expod}
f(x;\theta)=e^{-\theta x+\log\theta},\ \ \theta>0,
\end{equation}
The exponential distribution belongs to the single-parameter exponential family with $T(x)=-x$, $\psi(\theta)=-\log\theta$. The mean of exponential distribution is $\frac{1}{\theta}$. In this case, (\ref{thetav}) is simplified to
\begin{equation}\label{expoh}
\theta_1=\tau-\tau\frac{d}{
\sqrt{n}},\ \ \ \  \theta_2=\tau+\tau\frac{d}{\sqrt{n}}.
\end{equation}
Since $\frac{d^2\psi(\tau)}{d\theta^2}=\frac{1}{\tau^2}$, then $\boldsymbol\tau_{\boldsymbol d}^TD^2\psi(\boldsymbol\tau)\boldsymbol\tau_{\boldsymbol d}=d^2$. On the region $d^2<1$, there is no consistent. 

\noindent{\bf Example 2 (Normal distribution).}
The normal distribution   $N(\mu,\sigma^2)$ has density
\[f(x;\mu,\sigma^2)=\frac{1}{\sqrt{2\pi}}e^{-\frac{(x-\mu)^2}{2\sigma^2}}.\]
It is a member of the 2-parameter exponential family with $h(x)=\frac{1}{\sqrt{2\pi}}$, $T_1(x)=x$, $T_2(x)=x^2$, $\theta_1=\frac{\mu}{\sigma^2}$, $\theta_2=-\frac{1}{2\sigma^2}$ and $\psi(\theta_1,\theta_2)=-\frac{\theta_1^2}{4\theta_2}-\frac{1}{2}\log(-2\theta_2)$. Straightforward calculation yields
\[D^2\psi(\boldsymbol\theta)=
\begin{bmatrix}
-\frac{1}{2\theta_2} & \frac{\theta_1}{\theta_2^3}\\
\frac{\theta_1}{\theta_2^3} & -\frac{\theta_1^2}{2\theta_2^3}+\frac{1}{2\theta_2} 
\end{bmatrix}.
\]
Then 
\[\boldsymbol\tau_{\boldsymbol d}^TD^2\psi(\boldsymbol\tau)\boldsymbol\tau_{\boldsymbol d}=-\frac{1}{2\tau_2}\tau_1^2d_1^2+2\frac{\tau_1}{\tau_2^3}\tau_1\tau_2d_1d_2+\left(-\frac{\tau_1^2}{2\tau_2^3}+\frac{1}{2\tau_2}\right)\tau_2^2d_2^2.\]
For $\boldsymbol d$ within the ellipsoid $\boldsymbol\tau_{\boldsymbol d}^TD^2\psi(\boldsymbol\tau)\boldsymbol\tau_{\boldsymbol d}<1$, there is no consistent test.

 Theorem \ref{theorem:1} assumes the dimension $m$ of the exponential family to be finite. Actually, Theorem \ref{theorem:1} holds even for 
 infinite-dimensional exponential family defined in \cite{CS05}:
\[\mathcal{P}=\left\{p_f(x)=e^{f(x)-\psi(f)}q_0(x),\ x\in\Omega\subset \mathbb{R}: f\in\mathcal{F}\right\},\]
where $\mathcal{F}$ is a subspace of a reproducing kernel Hilbert space $\mathcal{H}$, $q_0(x)$ is a reference density function and
\[\mathcal{F}=\left\{f\in\mathcal{H}:e^{\psi(f)}<\infty\right\},\ \ \psi(f)=\log\int_{\Omega}e^{f(x)}q_0(x)dx.\]
The infinite-dimensional exponential family includes a very broad class of distributions (see \cite{CS05,SFGHK17}). 

Assume $P,Q\in\mathcal{P}$ are parametrized by $f_1$ and $f_2$ respectively. Let 
\[f_1=f-\frac{g}{\sqrt{n}},\hskip 1cm f_2=f+\frac{g}{\sqrt{n}},\]
where $g\in\mathcal{F}$ such that $f_1,f_2\in\mathcal{F}$. Then the hypotheses (\ref{hyp0}) can be reformulated as follows:
\begin{equation}\label{exphyp0}
    \textrm{$H_0:\|g\|_{\mathcal{H}}=0$ \,\,\,\,versus\,\,\,\, $H_1:\|g\|_{\mathcal{H}}\neq0$.}
\end{equation}
where $\|g\|_{\mathcal{H}}$ represents the norm of $g$ in the Hilbert space $\mathcal{H}$. In this case, the WSLMC test has the following form:
\begin{eqnarray*}
\mathcal{Z}_n=\frac{\sum_{(i_1,i_2,\dots,i_k)\in\mathcal{I}_n}\prod_{t=1}^k\frac{g(A_{i_ti_{t+1}})-\mathbb{E}_f(g(A_{i_ti_{t+1}}))}{\sqrt{n}}}{\sqrt{\frac{1}{2k}\left[D_F^2\psi(f)g^2\right]^k}},
\end{eqnarray*}
where $\mathbb{E}_f$ represents expectation with respect to density $p_f(x)\in\mathcal{P}$.
\begin{Theorem}\label{infexpo}
Suppose $D_F^5 \psi(f)$ exist and are uniformly bounded for all $f\in\mathcal{F}$. Then the following results hold.
\begin{itemize}
\item[(I)] If $D_F^2\psi(f)g^2<1$, any test is inconsistent. 
\item[(II)] If $D_F^2\psi(f)g^2>1$, the WSLMC test is consistent.
\end{itemize}
\end{Theorem}
To prove Theorem \ref{infexpo},
one needs to replace the partial derivatives by Fr$\Acute{e}$chet derivatives, and follows the line of proof of Theorem \ref{theorem:1}. In this sense, our result is pretty general.

\subsection{Information Loss of Dichotomizing Weighted Networks}\label{sec:2.2}

Existing community detection algorithms and statistical tests for community structure are mainly developed for binary edges (\cite{A17,AV14,AS17,ACB13,BM17,BS16,GL17a,L16,LR15,MS16,VA15,ZLZ12,ACB13}).  When a network is weighted, a common way is to convert the weighted network to a binary one. It has been empirically verified that there is information loss in the dichotomizing process (\cite{A14,AJC15,TB11}). However, it is unclear how much information is lost quantitatively.  In this subsection, we provide the first formal quantification of information loss in the context of hypothesis testing via statistical limits.

Given a fixed real number $t_0$, the weights $A_{ij}, i<j$ can be naturally dichotomized to binary ones $\tilde{A}_{ij}$ as follows
\[\tilde{A}_{ij}=I[A_{ij}>t_0],\ \ 1\leq i<j\leq n.\]
That is, all the weights smaller than $t_0$ are discarded and weights larger than $t_0$ are converted to ones (\cite{A14,AJC15,TB11}). In this way, the weighted work $A$ is converted to an unweighted network $\tilde{A}$. The network $\tilde{A}$ inherits the community structure of $A$. Given $\boldsymbol\sigma$ and $\boldsymbol\theta(\sigma_i,\sigma_j)=\boldsymbol\tau-\frac{\boldsymbol\tau_{\boldsymbol d}}{\sqrt{n}}\sigma_i\sigma_j$, the  probability of the presence of an edge in $\tilde{A}$ is
\[p_{ij}(\boldsymbol\sigma)=\mathbb{P}(\tilde{A}_{ij}=1|\boldsymbol\sigma)=\mathbb{P}(A_{ij}> t_0|\boldsymbol\sigma)=\int_{t_0}^{\infty}h(x)e^{\sum_{t=1}^m\theta_t(\sigma_i,\sigma_j)T_i(x)-\psi(\boldsymbol\theta(\sigma_i,\sigma_j))}dx.\]
Hence, the hypotheses are still the same as (\ref{exphyp0}).
Under $H_0$, there is no community structure and the edge presence probability is
\[p_0=\int_{t_0}^{\infty}h(x)e^{\boldsymbol\tau^T\boldsymbol T(x)-\psi(\boldsymbol\tau)}dx.\]

To get the sharp testing limit, we propose the signed-long-cycle (SLC) test statistic $\mathcal{R}_n$ as
\begin{eqnarray*}
\mathcal{R}_n=\frac{\sum_{(i_1,i_2,\dots,i_k)\in\mathcal{I}_n}\prod_{t=1}^k\left(\tilde{A}_{i_ti_{t+1}}-p_0\right)}{\sqrt{\frac{(k-1)!}{2}\binom{n}{k}\left[p_0(1-p_0)\right]^k}}, \hskip 1cm \  k=\log\log\log n.
\end{eqnarray*}
The SLC test rejects $H_0$ if $|\mathcal{R}_n|>C$ for some constant $C$ dependent on the type I error.

\begin{Theorem}\label{expoloss}
Suppose $\partial^{\alpha}\psi$, for $\alpha$ with $|\alpha|=5$, exist and are uniformly bounded for all $\boldsymbol\theta\in\Theta$. 
Let
\[\boldsymbol a(t_0)=\int_{t_0}^{\infty}h(x)e^{\boldsymbol\tau^T\boldsymbol T(x)-\psi(\boldsymbol\tau)}\left(-\boldsymbol T(x)+D\psi(\boldsymbol\tau)\right)dx.\]
Then the following results hold.
\begin{itemize}
\item[(I)] If $\frac{(\boldsymbol a(t_0)^T\boldsymbol\tau_{\boldsymbol d})^2}{p_0(1-p_0)}<1$, any test is inconsistent. 
\item[(II)] If $\frac{(\boldsymbol a(t_0)^T\boldsymbol\tau_{\boldsymbol d})^2}{p_0(1-p_0)}>1$, the SLC test is consistent.
\end{itemize}
\end{Theorem}

For the dichotomized network, the sharp limit for existence of consistent test is given by the quantity $\frac{(\boldsymbol a(t_0)^T\boldsymbol\tau_{\boldsymbol d})^2}{p_0(1-p_0)}$. 
In general, $\frac{(\boldsymbol a(t_0)^T\boldsymbol\tau_{\boldsymbol d})^2}{p_0(1-p_0)}$ differs from $\boldsymbol\tau_{\boldsymbol d}^TD^2\psi(\boldsymbol\tau)\boldsymbol\tau_{\boldsymbol d}$ in Theorem \ref{theorem:1}, since it depends on $t_0$. The difference between them can be considered as a measure of information loss incurred by dichotomizing a weighted network in the context of hypothesis testing. In this sense, our result provides the first theoretical characterization of loss of information.

In the general case, it is not immediately clearly what is the difference between $\frac{(\boldsymbol a(t_0)^T\boldsymbol\tau_{\boldsymbol d})^2}{p_0(1-p_0)}$ and  $\boldsymbol\tau_{\boldsymbol d}^TD^2\psi(\boldsymbol\tau)\boldsymbol\tau_{\boldsymbol d}$. For better illustration, we restrict $P,Q$ to be exponential distributions.
In this case, $p_0=e^{-\tau t_0}$ , $a(t_0)=t_0p_0$ and $\frac{(\boldsymbol a(t_0)^T\boldsymbol\tau_{\boldsymbol d})^2}{p_0(1-p_0)}=\frac{\tau^2t_0^2d^2}{e^{\tau t_0}-1}$. Then the following corollary follows.

\begin{Corollary}\label{inforloss}
Suppose $P,Q$ are exponential distributions given by (\ref{expod}), (\ref{expoh}) and $t_0>0$ is a fixed constant. Then the following results hold.
\begin{itemize}
\item[(I).] If $d^2<\frac{e^{\tau t_0}-1}{\tau^2t_0^2}$, any test is inconsistent. 
\item[(II).] If $d^2>\frac{e^{\tau t_0}-1}{\tau^2t_0^2}$, the SLC test is consistent.
\end{itemize}
\end{Corollary}

\begin{figure}
    \centering
    \includegraphics[width=8cm,height=7cm]{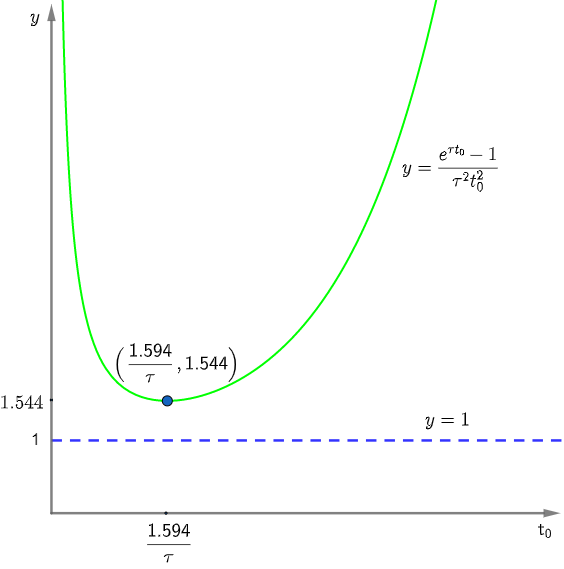}
    \caption{Comparison of the testing radius under exponential distributions. Blue dashed: testing radius ($\equiv 1$) under weighted network; Green solid: testing radius ($y$) under dichotomized weighted network. The latter is displayed as a function of threshold $t_0$ which achieves minimum at $t_0=\frac{1.594}{\tau}$. 
    }
    \label{lossFigure}
\end{figure}

For graphs with weights following the exponential distribution, the region without consistent test is $d^2<1$ by Example 1. For the dichotomized network, the region without consistent test is $d^2<\frac{e^{\tau t_0}-1}{\tau^2t_0^2}$ by Corollary \ref{inforloss}. Note that $\min_{x>0}\frac{e^{x}-1}{x^2}=1.544$ as shown in Figure \ref{lossFigure}. Hence, $\frac{e^{\tau t_0}-1}{\tau^2t_0^2}\geq1.544>1$ for any $\tau t_0>0$ and dichotomizing weighted exponential network to binary network always enlarges the region where no consistent test exists. This reflects the loss of information. For fixed $\tau$,  $t_0=\frac{1.594}{\tau}$ leads to the least loss of information. This interesting finding theoretically confirms the intuition that the dichotomy threshold should not be too small or too large. 

\subsection{A Stronger Result}\label{sec:2.3}

The result (I) in Theorem \ref{theorem:1} only guarantees the inconsistency of any statistical test, without providing further description of 
their performances. This can actually be made stronger under an additional ODE condition.
The following theorem says that, under such a condition, all tests perform not better than random guess. 
\begin{Theorem}\label{contiguity}
Suppose $\partial^{\alpha}\psi$, for $\alpha$ with $|\alpha|=5$, exist and are uniformly bounded for $\boldsymbol\theta\in\Theta$.
Moreover, $\psi$ satisfies the following differential equation
\begin{equation}\label{ODE}
\left(\sum_{|\alpha|=2}\frac{\partial^{\alpha}\psi}{\alpha!}(\boldsymbol\tau)\left(\boldsymbol\tau_{\boldsymbol d}\right)^{\alpha}\right)^2+3\sum_{|\alpha|=4}\frac{\partial^{\alpha}\psi}{\alpha!}(\boldsymbol\tau)\left(\boldsymbol\tau_{\boldsymbol d}\right)^{\alpha}=0.
\end{equation}
Then any test is not better than random guessing if $\boldsymbol\tau_{\boldsymbol d}D^2\psi(\boldsymbol\tau)\boldsymbol\tau_{\boldsymbol d}<1$.
\end{Theorem}

The following example provides distributions that satisfy condition (\ref{ODE}).

\noindent
{\bf Example 3}. For single-parameter exponential family, (\ref{ODE}) is reduced to
\[\left(\frac{d^2\psi}{d\eta^2}(\tau)\right)^2+\frac{1}{2}\frac{d^4\psi}{d\eta^4}(\tau)=0.\]
Consider the Gamma distribution with  density given by
\[f(x;\lambda,\eta)=\frac{\eta^{\lambda}}{\Gamma(\lambda)}x^{\lambda-1}e^{-\eta x}.\]
When $\lambda=3$, the Gamma distribution belong to the single-parameter exponential family with $\psi(\eta)=-3\log\eta$. Then direct calculation yields
\[\frac{d^2\psi}{d\eta^2}(\tau)=\frac{3}{\tau^2},\hskip 1cm \frac{d^4\psi}{d\eta^4}(\tau)=-\frac{18}{\tau^4},\]
which satisfies (\ref{ODE}). Similarly, the inverse Gamma distribution with scale parameter 2 satisfies (\ref{ODE}).

\section{Practical Tests}\label{practicalsec}

The WSLMC proposed in Section \ref{sec:2} is not practically feasible since it involves unknown quantities.
Specifically, it is unclear which exponential distribution the weights follow, and the vector $\boldsymbol T(x)$ is unknown as well. However, the  WSLMC test still provides useful insights about how to construct a powerful test. Note that for exponential distribution, $\boldsymbol T(x)=-x$; and for the normal distribution, $\boldsymbol T(x)=(x,x^2)^T$ . This hints us that we should consider higher moments of weights when proposing novel test statistics. Based on this observation,  we propose more practical tests that work for a broad class of distributions $P$ and $Q$.

To be convenient, we reformulate the hypotheses (\ref{hyp0}) as follows. Given a fixed positive integer $m$ and $\Theta\subset \mathbb{R}^m$, let $f(x;\boldsymbol\mu)$ be a probability density with parameter $\boldsymbol\mu=(\mu_1,\mu_2,\dots,\mu_m)^T\in\Theta$, where $\mu_t$ ($t\in\{1,2,\dots,m\}$) is the $t$th moment. Define the $m$-parameter distribution family as
\[\mathcal{F}_m=\{f(x;\boldsymbol\mu)|\boldsymbol\mu=(\mu_1,\mu_2,\dots,\mu_m)^T\in\Theta\subset \mathbb{R}^m\}.\]
 Suppose $P$ and $Q$ belong to $\mathcal{F}_m$ and has parameter $\boldsymbol\mu_1$ and $\boldsymbol\mu_2$ respectively. Let

\[
\boldsymbol\mu_1=\boldsymbol\mu+\frac{\boldsymbol\mu_{\boldsymbol d}}{\sqrt{n}},\,\,\,\, \boldsymbol\mu_2=\boldsymbol\mu-\frac{\boldsymbol\mu_{\boldsymbol d}}{\sqrt{n}},
\]
where $\boldsymbol\mu=(\mu_1,\mu_2,\dots,\mu_m)^T$, $\boldsymbol d=(d_1,d_2,\dots,d_m)^T$ and $\boldsymbol\mu_{\boldsymbol d}=(\mu_1d_1,\mu_2d_2,\dots,\mu_md_m)^T$. Then (\ref{hyp0}) can be rewritten as
\begin{equation}\label{exphyp}
    \textrm{$H_0:\|\boldsymbol{d}\|_2=0$\,\,\,\, versus\,\,\,\, $H_1:\|\boldsymbol{d}\|_2\neq 0$.}
\end{equation}

Let $\boldsymbol M(x)=(x,x^2,\dots,x^m)^T$ and $k\ge3$ be a positive integer. 
Define \[\overline{A}_t=\frac{1}{\binom{n}{2}}\sum_{1\leq i<j\leq n}A_{ij}^t,\,\,\,\, t=1,2,\dots,m,\] 
and $\overline{\boldsymbol M}( \boldsymbol A)=(\overline{A}_1,\overline{A}_2,\dots,\overline{A}_m)^T$. The sample covariance of $\boldsymbol M(A_{12})$ is
\[\boldsymbol S^2=\frac{1}{\binom{n}{2}}\sum_{1\leq i<j\leq n}(\boldsymbol M(A_{ij})-\overline{\boldsymbol M}( \boldsymbol A))(\boldsymbol M(A_{ij})-\overline{\boldsymbol M}( \boldsymbol A))^T.\]

We propose several tests based on the number of cycles. The first test is based on the signed-long-mixture-cycle (SLMC) test statistic defined as
\begin{eqnarray*}
\mathcal{T}_n=\frac{\sum_{(i_1,i_2,\dots,i_k)\in\mathcal{I}_n}\prod_{t=1}^k\boldsymbol 1^T\left(\boldsymbol M(A_{i_ti_{t+1}})-\overline{\boldsymbol M}( \boldsymbol A)\right)}{\sqrt{\frac{(k-1)!}{2}\binom{n}{k}\big(\boldsymbol 1^T\boldsymbol S^21\big)^k}}.
\end{eqnarray*}

\begin{Theorem}\label{prac1}
Suppose $3\leq k=O\left(\log\log\log n\right)$ and all the moments of $P$and $Q$ exist. Then under $H_0$, $\mathcal{T}_n$ converges in distribution to the standard normal distribution as $n$ goes to infinity.
\end{Theorem}

Based on Theorem \ref{prac1},
the SLMC test rejects $H_0$ if $|\mathcal{T}_n|>Z_{\frac{\gamma}{2}}$, where $Z_{\frac{\gamma}{2}}$ is the $100\frac{\gamma}{2}\%$ quantile of the standard normal distribution. Since the number of cycles with length $k$ can be expressed as a function of the trace of $A^k$ and the number of $k$ walks, the computation complexity of $\mathcal{T}_n$ is at most $O(n^3k)$. Hence, $\mathcal{T}_n$ is a practical test statistic. We point out that the condition all the moments of $P,Q$ exist can be relaxed to that the $4m$-th moments exist. The current proof of Theorem \ref{prac1} employs the method of moment which requires all the moments are finite. An alternative proof is to use the Martingale central limit theorem, which only requires finite $4m$-th moment.

\begin{Theorem}\label{prac2}
Suppose $3\leq k=O\left(\log\log\log n\right)$ , $2m$-th moments of $P$ and $Q$ exist and $\max_{1\leq t\leq m}{d_t}=o(\sqrt{n})$. Then under $H_1$, $\mathcal{T}_n=\frac{1}{\sqrt{2k}}\big(\frac{\boldsymbol 1^T\boldsymbol\mu_{\boldsymbol d}}{\sqrt{\boldsymbol 1^T\boldsymbol \Sigma \boldsymbol 1}}\big)^k\left(1+o_p(1)\right)$.
\end{Theorem}

Based on Theorem \ref{prac2}, the power of the SLMC test approaches one as $n\rightarrow\infty$ if $\frac{1}{\sqrt{2k}}\big(\frac{\boldsymbol 1^T\boldsymbol\mu_{\boldsymbol d}}{\sqrt{\boldsymbol 1^T\boldsymbol \Sigma \boldsymbol 1}}\big)^k\rightarrow\infty$. In this case, the SLMC test with larger $k$ may achieve higher power. If $\max_{1\leq t\leq m}{d_t}\rightarrow\infty$, the power can tend to one for finite $k$. Most importantly, even when $\max_{1\leq t\leq m}{d_t}$ is bounded, our test can still have asymptotic power one whenever $k\rightarrow\infty$ and $k=O\left(\log\log\log n\right)$. In this sense, the  SLMC test is almost optimal. The condition $\max_{1\leq t\leq m}{d_t}=o(\sqrt{n})$ is just to simplify the order of $\mathcal{T}_n$ under $H_1$. When $\max_{1\leq t\leq m}{d_t}=c\sqrt{n}$ for some constant $c>0$, the order of $\mathcal{T}_n$ has a tedious expression but the power still converges to one. Theorem \ref{prac2} only requires finite $2m$-th moments of $P$ and $Q$, since we did not pursue the asymptotic distribution of $\mathcal{T}_n$. It suffices to get the order of $\mathcal{T}_n$ under $H_1$ for power analysis.

The SLMC test statistic employs all the $m$ moments. Alternatively,
we can also use a single moment to construct a test statistic. Let  $l\in\{1,2,\dots,m\}$ and
\[ S_l^2=\frac{1}{\binom{n}{2}}\sum_{i<j}(A_{ij}^l-\overline{A}_l)^2.\]
Define the signed long-cycle (SLC) test statistic as
\begin{eqnarray*}
\mathcal{T}_{n,l}=\frac{\sum_{(i_1,i_2,\dots,i_k)\in\mathcal{I}_n}\prod_{t=1}^k\left(A_{i_ti_{t+1}}^l-\overline{A}_l\right)}{\sqrt{\frac{(k-1)!}{2}\binom{n}{k}S_l^{2k}}}.
\end{eqnarray*}

The following results follow by a similar proof of Theorem \ref{prac1} and Theorem \ref{prac2}.

\begin{Corollary}
Suppose $l\in\{1,2,\dots,m\}$ is a fixed integer, $3\leq k=O\left(\log\log\log n\right)$ and all the moments of $P$and $Q$ exist. Then under $H_0$, $\mathcal{T}_{n,l}$ converges in distribution to the standard normal distribution as $n$ goes to infinity.
Under $H_1$, if $d_l=o(\sqrt{n})$, then $\mathcal{T}_{n,l}=\frac{1}{\sqrt{2k}}\big(\frac{\mu_ld_l}{\sigma_l}\big)^k\left(1+o_p(1)\right)$. Here $\sigma_l^2=Var(A_{12}^l)$ under $H_0$.
\end{Corollary}

The SLC test rejects $H_0$ if $|\mathcal{T}_{n,l}|>Z_{\frac{\gamma}{2}}$, where $Z_{\frac{\gamma}{2}}$ is the $100\frac{\gamma}{2}\%$ quantile of the standard normal distribution.  
The power approaches one as $n\rightarrow\infty$ if $ \frac{1}{\sqrt{2k}}\big(\frac{\mu_ld_l}{\sigma_l}\big)^k\rightarrow\infty$. When $\boldsymbol\mu_{\boldsymbol d}$ only has a single nonzero component, the SCL test may have higher power than the SLMC test.

\section{Simulation and Application}\label{sec:4}

\subsection{Simulation}

In this subsection, we illustrate the performance of the
proposed tests and compare them with the spectral test (\cite{T18}) in various simulations. The nominal type I error is set to be 0.05. The empirical type I errors and powers are calculated based on 500 repetitions.

For each $l\in\{1,2,\dots,m\}$, the spectral test statistics  are defined as
\[\Lambda_{n,l}=n^{\frac{2}{3}}(\lambda_{l,max}-2),\,\,\,\,\Lambda_{n}=n^{\frac{2}{3}}(\lambda_{n}-2),\]
where
\[\lambda_{l,max}=\lambda_{max}\left(\frac{A^l-\overline{A}_l}{\sqrt{n}S_l}\right),\hskip 1cm
\lambda_n=\lambda_{max}\left(\frac{\boldsymbol M(A)-\overline{\boldsymbol M}( \boldsymbol A)}{\sqrt{n\boldsymbol 1^T\boldsymbol S^2\boldsymbol1}}\right).\]
Note that $\Lambda_{n,1}$ is given in \cite{T18}. The rejection threshold is calculated based on the Tracy-Widom law with degree of freedom one. 

We assess the performance of the cycle-based tests $\mathcal{T}_n$, $\mathcal{T}_{n,1}$, $\mathcal{T}_{n,2}$ and the spectral tests $\Lambda_{n,1},\Lambda_{n,2}$. The results for $\Lambda_n$ are not reported due to large type I errors. The random labels $\sigma_i$ $(1\leq i\leq n)$ are generated from Bernoulli distribution with success probability 0.5. Given $\boldsymbol\sigma=(\sigma_1,\sigma_2,\dots,\sigma_n)$, let
$\boldsymbol\lambda=(\lambda_1,\lambda_2)$ and $\boldsymbol\epsilon=(\epsilon_1,\epsilon_2)$, and
\begin{equation}\label{simu}
\mu_1=\lambda_1+\epsilon_1\sigma_i\sigma_j,\hskip 1cm \mu_2=\lambda_2+\epsilon_2\sigma_i\sigma_j,
\end{equation}

\begin{figure}
    \centering
    \includegraphics[width=15cm,height=6cm]{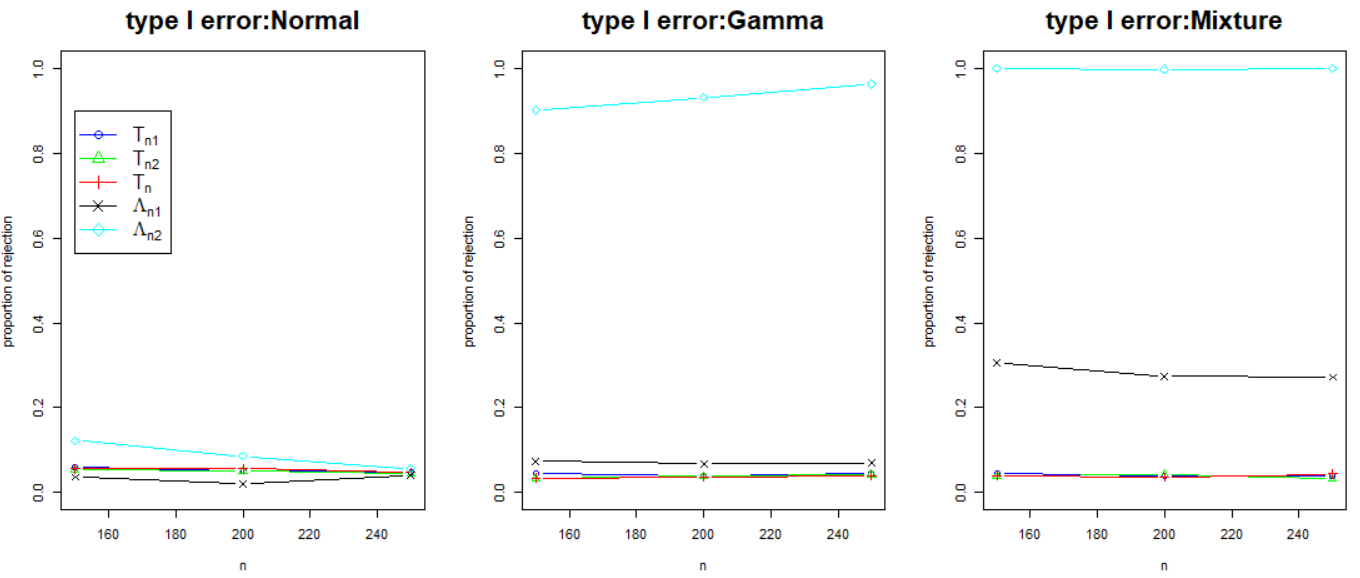}
\caption{Type I errors of the spectral tests and cycle-based tests in three models.}
    \label{fig:size}
\end{figure}

\begin{figure}
    \centering
    \includegraphics[width=15cm,height=6cm]{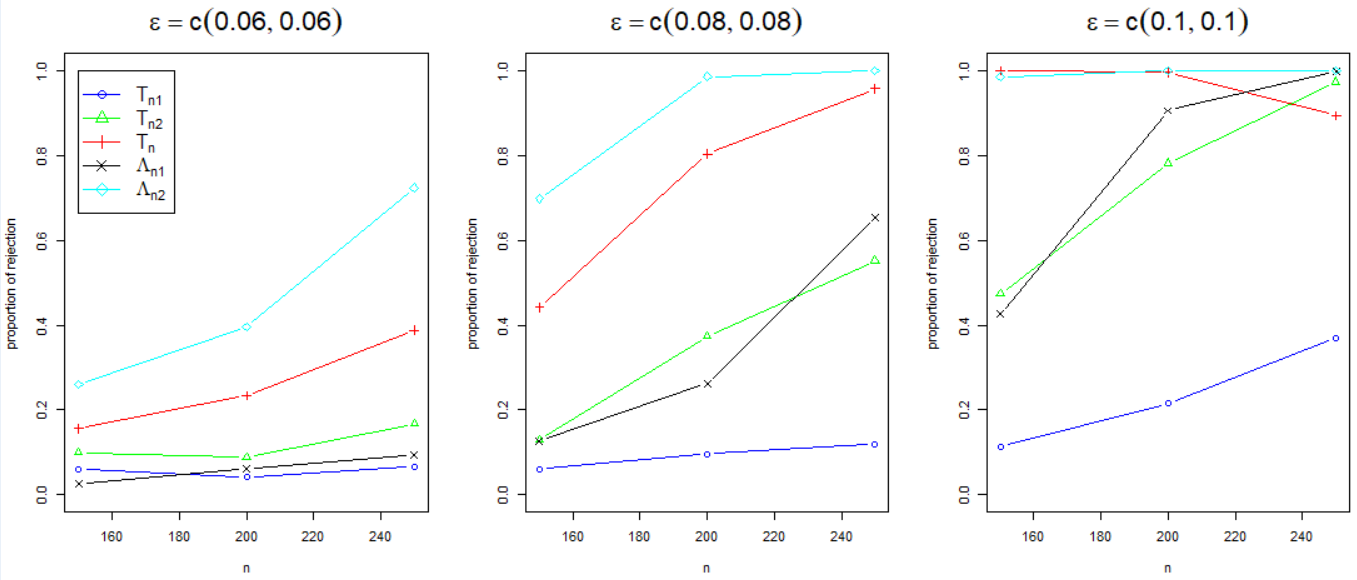}
    \includegraphics[width=15cm,height=6cm]{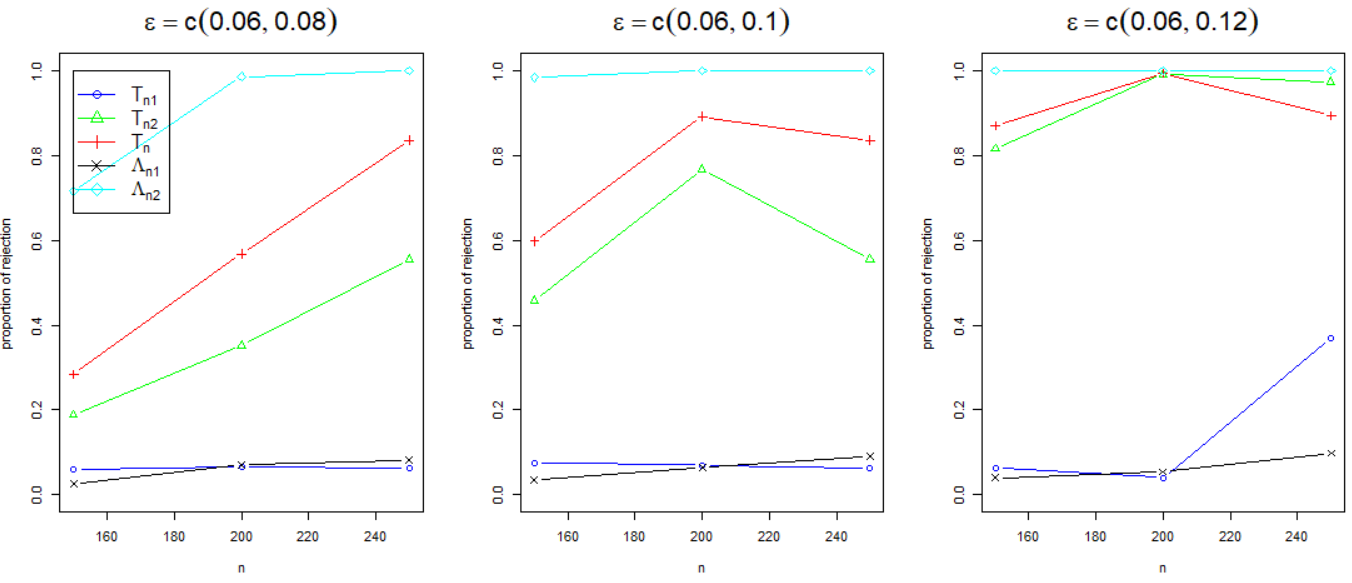}
  \includegraphics[width=15cm,height=6cm]{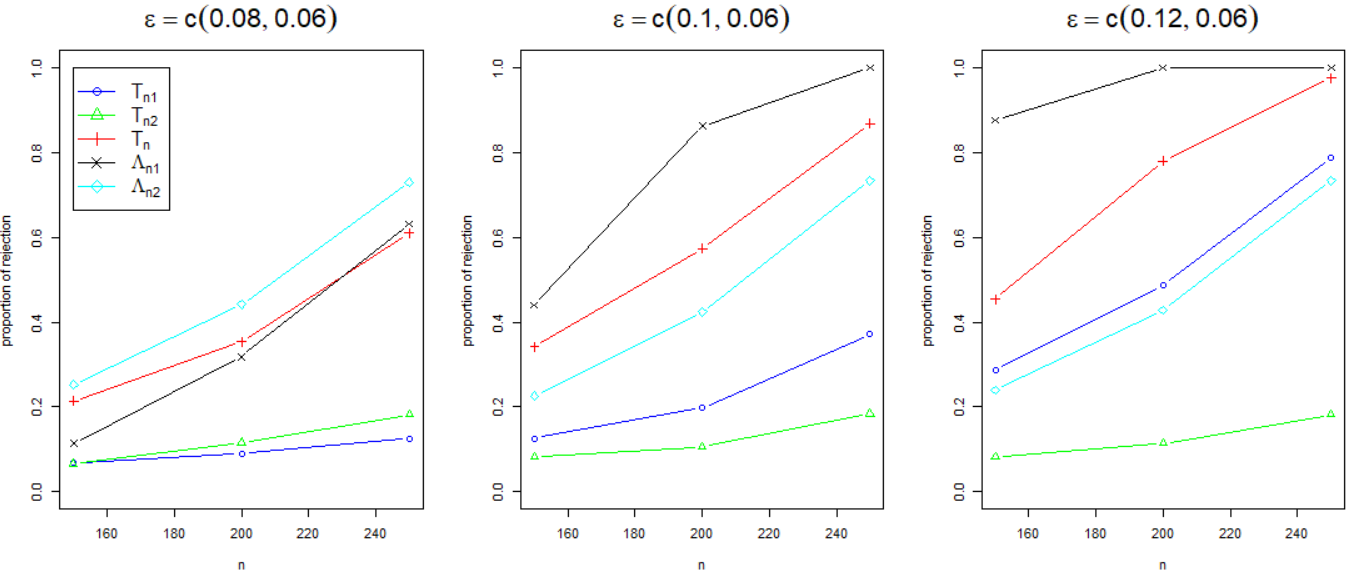}
\caption{Powers of the spectral tests and cycle-based tests when weights follow Normal distribution with $\lambda=(0,1)$.}
    \label{fig:normal}
\end{figure}

In the first simulation, we generate the weights from the normal distribution with mean $\mu_1$ and second moment $\mu_2$ defined in (\ref{simu}). Figure \ref{fig:size} summarizes the empirical type I errors and Figure \ref{fig:normal} visualizes the powers.

\begin{figure}
    \centering
    \includegraphics[width=15cm,height=6cm]{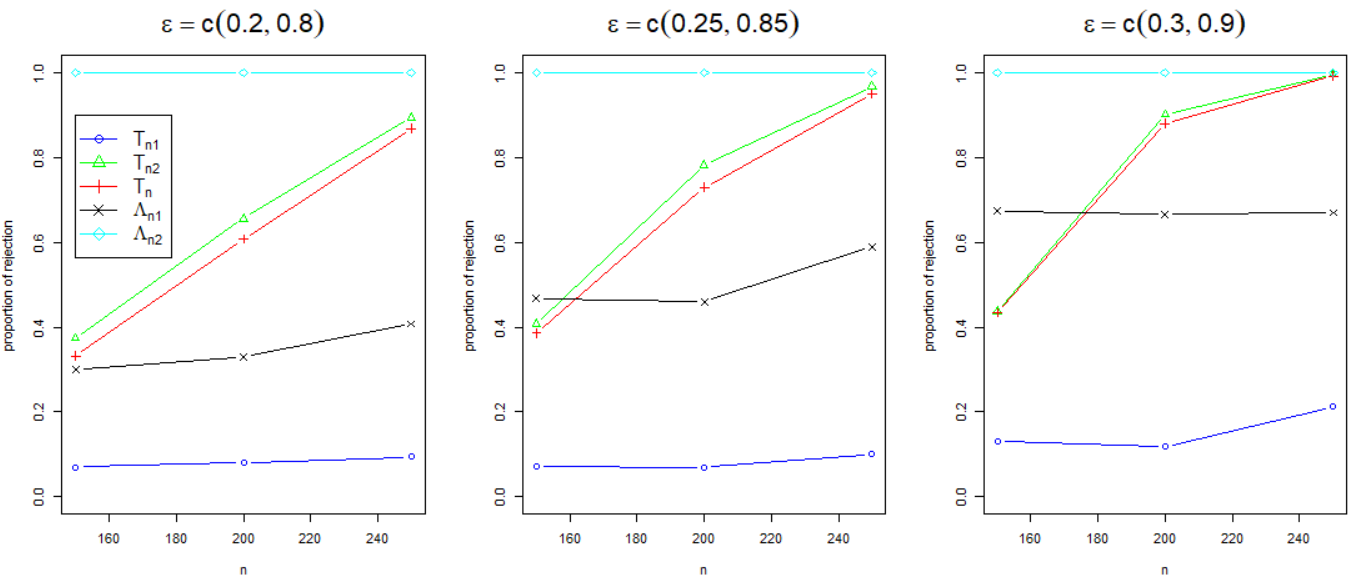}
    \includegraphics[width=15cm,height=6cm]{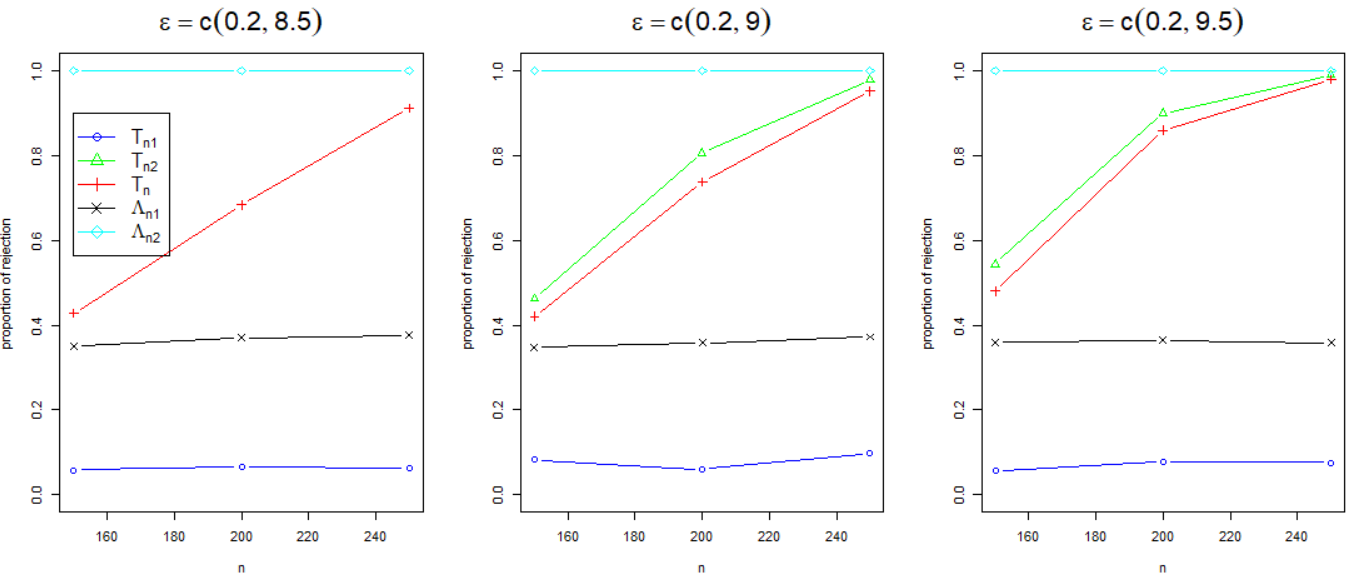}
  \includegraphics[width=15cm,height=6cm]{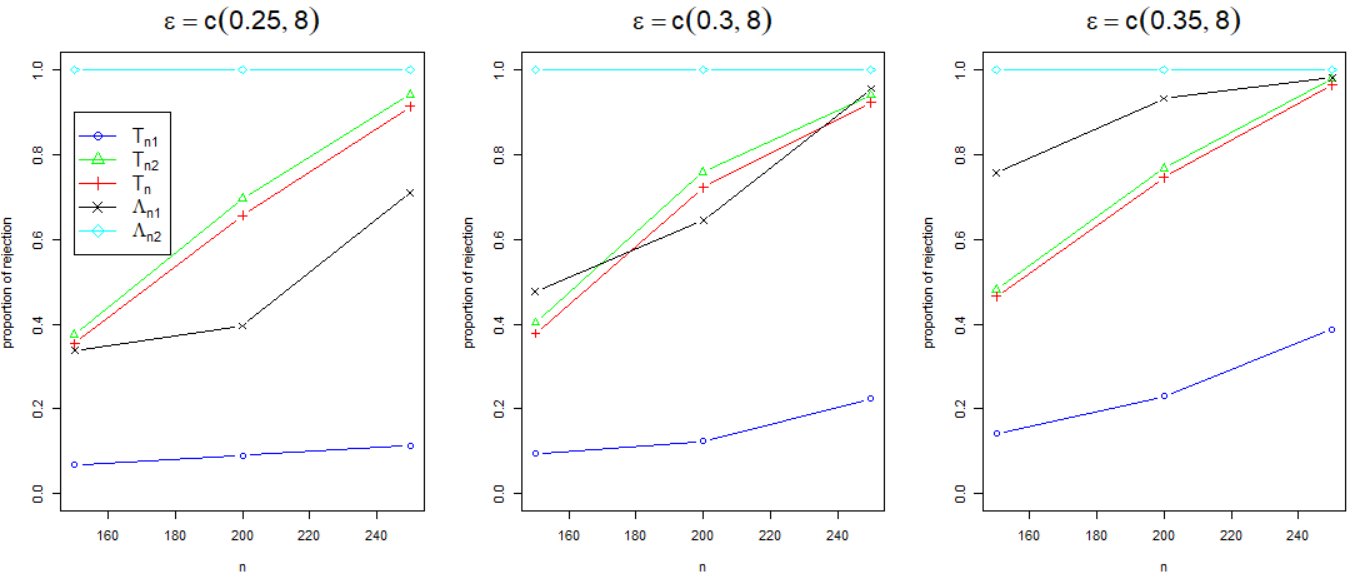}
\caption{Powers of the spectral tests and cycle-based tests when weights follow Gamma distribution with $\lambda=(4,28)$.}
    \label{fig:gamma}
\end{figure}

In the second simulation, the weights are assumed to follow the Gamma distribution with mean and second moment given in (\ref{simu}). The density of Gamma distribution is given by
\[f(x;\lambda,\theta)=\frac{1}{\Gamma(\lambda)}x^{\lambda-1}e^{-\frac{x}{\theta}},\]
where $\lambda,\theta$ are functions of $\mu_1,\mu_2$:
\[\lambda=\frac{\mu_1^2}{\mu_2-\mu_1^2},\hskip 1cm \theta=\frac{\mu_2-\mu_1^2}{\mu_1}.\]  Figure \ref{fig:size} presents the empirical type I errors and Figure \ref{fig:gamma} plots the powers.

\begin{figure}
    \centering
    \includegraphics[width=15cm,height=6cm]{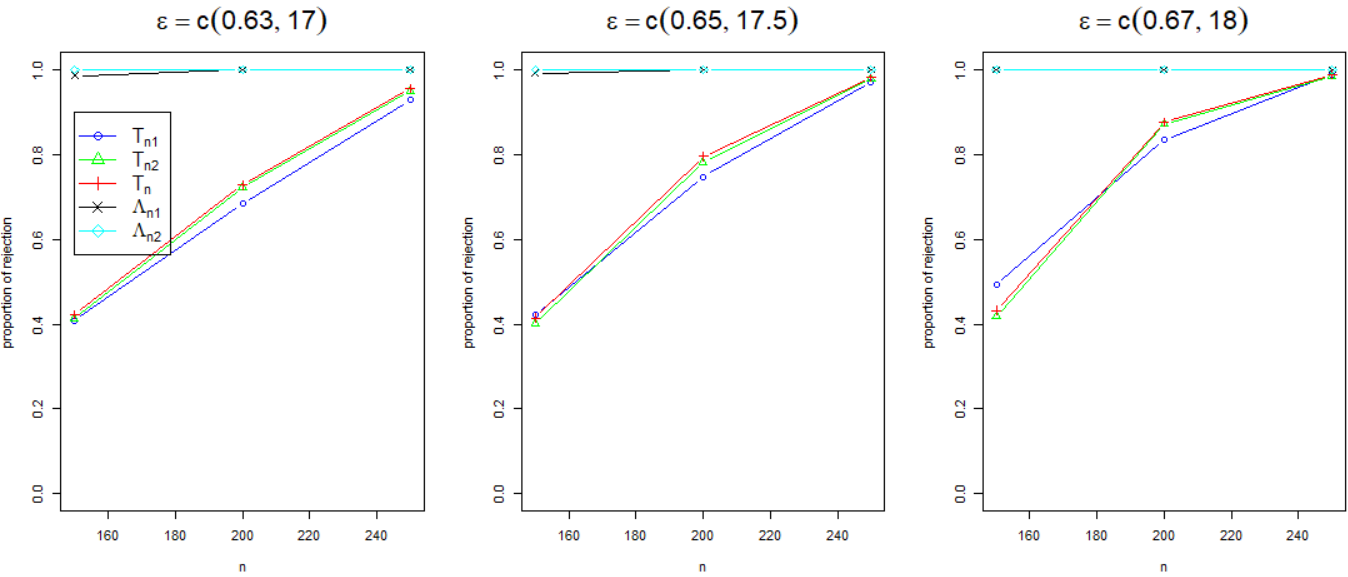}
    \includegraphics[width=15cm,height=6cm]{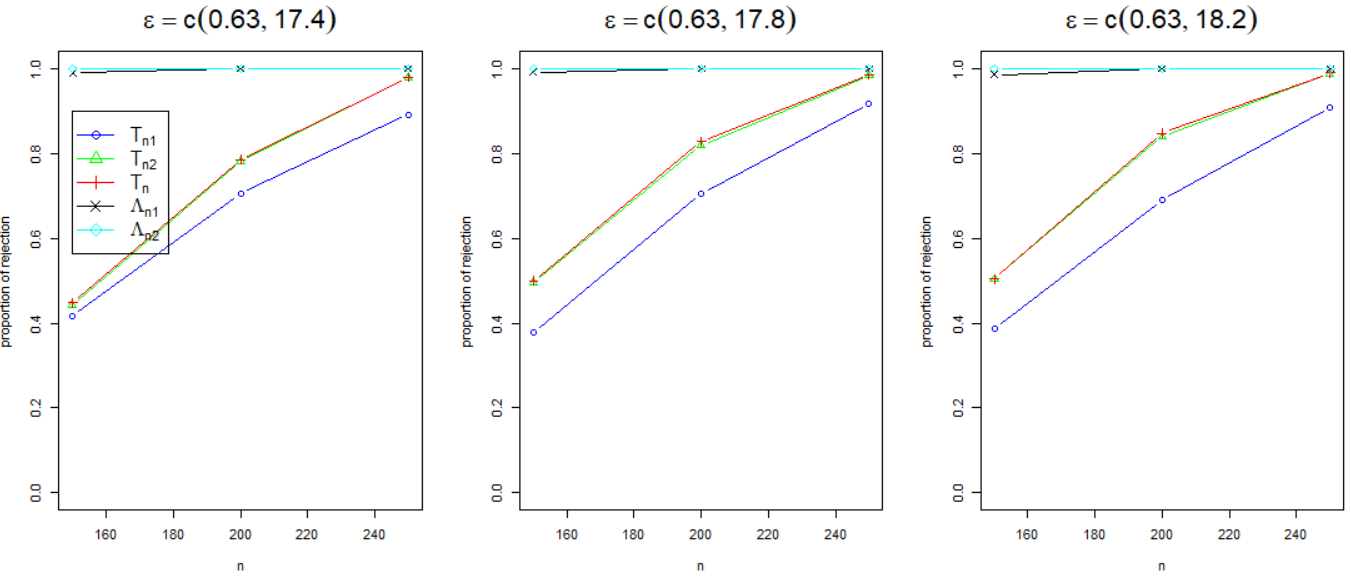}
  \includegraphics[width=15cm,height=6cm]{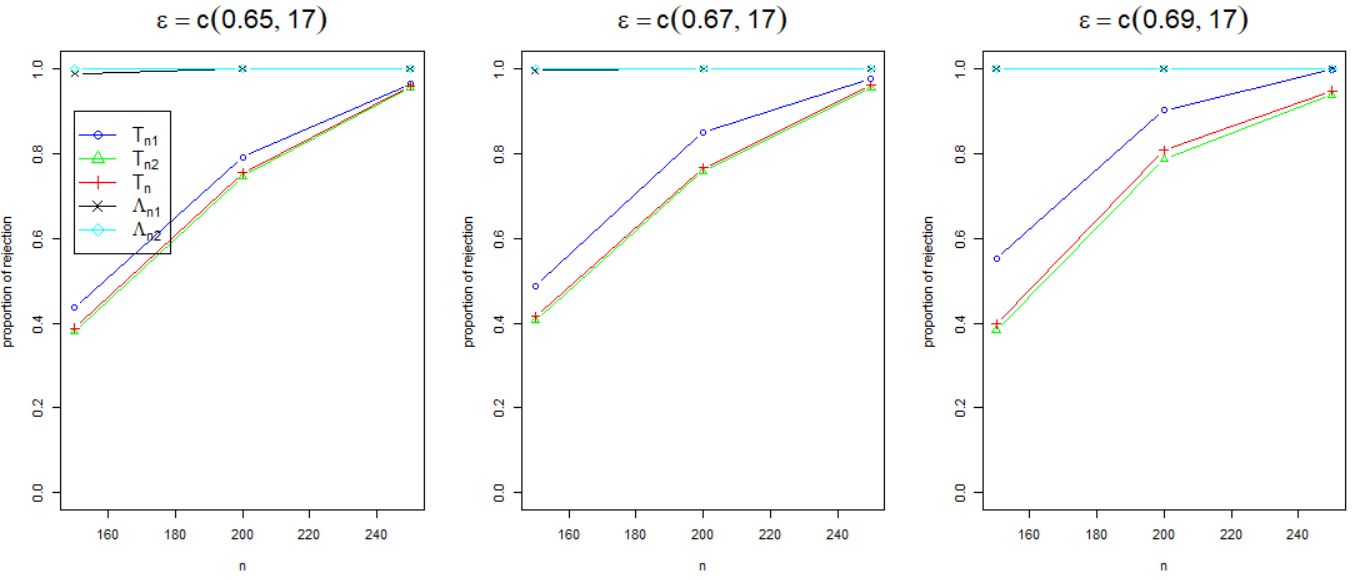}
\caption{Powers of the spectral tests and cycle-based tests when weights follow mixture exponential distribution with $\lambda=(3.6, 36)$.}
    \label{fig:mixture}
\end{figure}

In the last simulation, we generate weights from the 
mixture of two exponential distributions with density given by
\[f(x;\lambda_1,\lambda_2)=0.5\lambda_1e^{-\lambda_1x}+0.5\lambda_2e^{-\lambda_2x},\]
where
\[\lambda_1=\frac{2}{2\mu_1+\sqrt{2\mu_2-4\mu_1^2}},\hskip 1cm \lambda_2=\frac{2}{2\mu_1-\sqrt{2\mu_2-4\mu_1^2}}.\]
Note that this distribution does not belong to exponential family.
  The empirical type I errors are plotted in Figure \ref{fig:size} and Figure \ref{fig:mixture} visualizes the powers.

This simulation study has the following indications: $(I)$ The cycle-based tests with $k=3$ converge much faster than the spectral tests for non-Gaussian weights, since the cycle-based tests have empirical type I errors close to the nominal for smaller $n$. $(II)$ The empirical powers increase as the differences ($\epsilon_1$ or $\epsilon_2$) of parameters get larger. $(III)$ No single test can dominate the others in all cases. In summary, this simulation highlights the necessity to incorporate higher order moments of weights in test statistic and the proposed tests have significant advantage over spectral tests for smaller $n$.

\subsection{Animal Social Network Data}

In this subsection, we apply the proposed tests to an animal social network ``aves-wildbird-network-5'' publicly available at \cite{RA15}. In this network, each bird is a node and edges are weighted by interaction between two birds. This network has 145 nodes and 2512 weighted edges. We apply test statistics $\mathcal{T}_n$, $\mathcal{T}_{n,1}$, $\mathcal{T}_{n,2}$ with $k=3$ to this network and the calculated test statistics are listed in Table \ref{realdata}. With type I error $0.05$, all these three tests reject the null hypothesis, which indicates the network contains community structure.

\begin{table}[h]
\centering
\caption{Calculated test statistics.}
\centering
\begin{tabular}{ |p{1.5cm}|p{1.5cm}|p{1.5cm}|}
\hline
 $\mathcal{T}_{n,1}$ & $\mathcal{T}_{n,2}$  &  $\mathcal{T}_{n}$ \\
\hline
    225.7204  &  47.9405  &  174.9434     \\ 
\hline
\end{tabular}
\label{realdata}
\end{table}

\section{Proof of main result}\label{sec:app}

\subsection{Proof of Theorem \ref{theorem:1}}

 For the proof of part (I) of Theorem \ref{theorem:1}, we will use the second moment methods. Specifically, we show that the second moment of the likelihood ratio under $H_0$ is bounded if $2\sum_{|\alpha|=2}\frac{\partial^{\alpha}\psi}{\alpha!}(\boldsymbol\tau)\left(\boldsymbol\tau_{\boldsymbol d}\right)^{\alpha}<1$. For the proof of part (II) of Theorem \ref{theorem:1}, we prove the WSLMC test has asymptotic power one if $2\sum_{|\alpha|=2}\frac{\partial^{\alpha}\psi}{\alpha!}(\boldsymbol\tau)\left(\boldsymbol\tau_{\boldsymbol d}\right)^{\alpha}>1$. For convenience, we will let $m=2$. The proof for general $m$ is exactly the same as $m=2$.

\noindent
{\bf Proof of Theorem \ref{theorem:1} (I):} The proof strategy is to show the second moment of the likelihood ratio under $H_0$ is bounded. Given random label vector $\boldsymbol\sigma=(\sigma_1,\sigma_2,\dots,\sigma_n)$, the parameters of the distribution of $A_{ij}$ can be concisely written as
\begin{equation}\label{theta1}
\theta_{t,ij}(\boldsymbol\sigma)=\tau_{t}-\tau_{t}\frac{d_t}{\sqrt{n}}\sigma_i\sigma_j,\hskip 1cm t=1,2.
\end{equation}
Let $\boldsymbol\theta_{ij}(\boldsymbol\sigma)=(\theta_{1,ij}(\boldsymbol\sigma),\theta_{2,ij}(\boldsymbol\sigma))$ and   $\boldsymbol\tau=(\tau_{1},\tau_{2})$.
Then the likelihood ratio $L_n$ is equal to
\begin{eqnarray*}
L_n=\frac{1}{2^n}\sum_{\boldsymbol\sigma\in\{\pm\}^n}\left[\prod_{1\leq i<j\leq n}\frac{h(A_{ij})\exp\{\theta_{1,ij}(\boldsymbol\sigma) T_1(A_{ij})+\theta_{2,ij}(\boldsymbol\sigma) T_2(A_{ij})-\psi(\boldsymbol\theta_{ij}(\boldsymbol\sigma))\}}{h(A_{ij})\exp\{\tau_1 T_1(A_{ij}))+\tau_2 T_2(A_{ij}))-\psi(\boldsymbol\tau)\}}\right].
\end{eqnarray*}
Let $\boldsymbol\eta$ be an independent copy of $\boldsymbol\sigma$. Then under $H_0$, the second moment of $L_n$ is equal to
\begin{eqnarray}\nonumber
&&\mathbb{E}[L_n^2]\\ \nonumber
&=&\frac{1}{4^n}\sum_{\boldsymbol\sigma,\boldsymbol\eta\in\{\pm\}^n}\mathbb{E}\Bigg[\prod_{1\leq i<j\leq n}\frac{h(A_{ij})\exp\{\theta_{1,ij}(\boldsymbol\sigma) T_1(A_{ij})+\theta_{2,ij}(\boldsymbol\sigma) T_2(A_{ij})-\psi(\boldsymbol\theta_{ij}(\boldsymbol\sigma))\}}{h(A_{ij})\exp\{\tau_1 T_1(A_{ij}))+\tau_2 T_2(A_{ij}))-\psi(\boldsymbol\tau)\}}\\ \nonumber
&&\times\prod_{1\leq i<j\leq n}\frac{h(A_{ij})\exp\{\theta_{1,ij}(\boldsymbol\eta) T_1(A_{ij})+\theta_{2,ij}(\boldsymbol\eta) T_2(A_{ij})-\psi(\boldsymbol\theta_{ij}(\boldsymbol\eta))\}}{h(A_{ij})\exp\{\tau_1 T_1(A_{ij}))+\tau_2 T_2(A_{ij}))-\psi(\boldsymbol\tau)\}}\Bigg]\\  \nonumber
&=&\mathbb{E}_{\boldsymbol\sigma,\boldsymbol\eta}\Bigg[\prod_{1\leq i<j\leq n}\int\frac{h(A_{ij})\exp\{\theta_{1,ij}(\boldsymbol\sigma) T_1(A_{ij})+\theta_{2,ij}(\boldsymbol\sigma) T_2(A_{ij})-\psi(\boldsymbol\theta_{ij}(\boldsymbol\sigma))\}}{\exp\{\tau_1 T_1(A_{ij}))+\tau_2 T_2(A_{ij}))-\psi(\boldsymbol\tau)\}}\\ \nonumber
&&\times \exp\{\theta_{1,ij}(\boldsymbol\eta) T_1(A_{ij})+\theta_{2,ij}(\boldsymbol\eta) T_2(A_{ij})-\psi(\boldsymbol\theta_{ij}(\boldsymbol\eta))\}dA_{ij}\Bigg]\\  \nonumber
&=&\mathbb{E}_{\boldsymbol\sigma,\boldsymbol\eta}\Bigg[\prod_{1\leq i<j\leq n}\int h(A_{ij})\exp\Big\{[\theta_{1,ij}(\boldsymbol\sigma)+\theta_{1,ij}(\boldsymbol\eta)-\tau_1] T_1(A_{ij})+[\theta_{2,ij}(\boldsymbol\sigma)+\theta_{2,ij}(\boldsymbol\eta)-\tau_2] T_2(A_{ij})\\ \nonumber
&&-\psi(\boldsymbol\theta_{ij}(\boldsymbol\sigma))-\psi(\boldsymbol\theta_{ij}(\boldsymbol\eta))+\psi(\boldsymbol\tau)\Big\}dA_{ij}\Bigg]\\  \label{exp2}
&=&\mathbb{E}_{\boldsymbol\sigma,\boldsymbol\eta}\left[\prod_{1\leq i<j\leq n}\exp\Big\{\psi\Big(\boldsymbol\theta_{ij}(\boldsymbol\sigma)+\boldsymbol\theta_{ij}(\boldsymbol\eta)-\boldsymbol\tau\Big) -\psi(\boldsymbol\theta_{ij}(\boldsymbol\sigma))-\psi(\boldsymbol\theta_{ij}(\boldsymbol\eta))+\psi(\boldsymbol\tau)\Big\}\right]
\end{eqnarray}
For $\sigma_i\sigma_j=1$ and $\eta_i\eta_j=1$, we have
\begin{eqnarray}\nonumber
&&\psi\Big(\boldsymbol\theta_{ij}(\boldsymbol\sigma)+\boldsymbol\theta_{ij}(\boldsymbol\eta)-\boldsymbol\tau\Big) -\psi(\boldsymbol\theta_{ij}(\boldsymbol\sigma))-\psi(\boldsymbol\theta_{ij}(\boldsymbol\eta))+\psi(\boldsymbol\tau)\\ \label{exp3}
&=&\psi\left(\tau_1-\frac{2\tau_1d_1}{\sqrt{n}},\tau_2-\frac{2\tau_2d_2}{\sqrt{n}}\right)-2\psi\left(\tau_1-\frac{\tau_1d_1}{\sqrt{n}},\tau_2-\frac{\tau_2d_2}{\sqrt{n}}\right)+\psi(\tau_1,\tau_2).
\end{eqnarray}
By Taylor expansion, we have
\begin{eqnarray}\nonumber
&&\psi\left(\tau_1-\frac{2\tau_1d_1}{\sqrt{n}},\tau_2-\frac{2\tau_2d_2}{\sqrt{n}}\right)\\ \nonumber
&=&\psi(\tau_1,\tau_2)-\frac{2\boldsymbol\tau_{\boldsymbol d}}{\sqrt{n}}D\psi(\boldsymbol\tau)+\sum_{|\alpha|=2}\frac{\partial^{\alpha}\psi}{\alpha!}(\boldsymbol\tau)\left(\frac{-2\boldsymbol\tau_{\boldsymbol d}}{\sqrt{n}}\right)^{\alpha}+\sum_{|\alpha|=3}\frac{\partial^{\alpha}\psi}{\alpha!}(\boldsymbol\tau)\left(\frac{-2\boldsymbol\tau_{\boldsymbol d}}{\sqrt{n}}\right)^{\alpha}\\ \label{exp4}
&&+\sum_{|\alpha|=4}\frac{\partial^{\alpha}\psi}{\alpha!}(\boldsymbol\tau)\left(\frac{-2\boldsymbol\tau_{\boldsymbol d}}{\sqrt{n}}\right)^{\alpha}+O\left(\frac{1}{n^2\sqrt{n}}\right).
\end{eqnarray}

\begin{eqnarray}\nonumber
&&\psi\left(\tau_1-\frac{\tau_1d_1}{\sqrt{n}},\tau_2-\frac{\tau_2d_2}{\sqrt{n}}\right)\\ \nonumber
&=&\psi(\tau_1,\tau_2)-\frac{\boldsymbol\tau_{\boldsymbol d}}{\sqrt{n}}D\psi(\boldsymbol\tau)+\sum_{|\alpha|=2}\frac{\partial^{\alpha}\psi}{\alpha!}(\boldsymbol\tau)\left(\frac{-\boldsymbol\tau_{\boldsymbol d}}{\sqrt{n}}\right)^{\alpha}+\sum_{|\alpha|=3}\frac{\partial^{\alpha}\psi}{\alpha!}(\boldsymbol\tau)\left(\frac{-\boldsymbol\tau_{\boldsymbol d}}{\sqrt{n}}\right)^{\alpha}\\ \label{exp5}
&&+\sum_{|\alpha|=4}\frac{\partial^{\alpha}\psi}{\alpha!}(\boldsymbol\tau)\left(\frac{-\boldsymbol\tau_{\boldsymbol d}}{\sqrt{n}}\right)^{\alpha}+O\left(\frac{1}{n^2\sqrt{n}}\right).
\end{eqnarray}
Hence, by (\ref{exp3}), (\ref{exp4}),(\ref{exp5}), we get
\begin{eqnarray}\nonumber
&&\psi\Big(\boldsymbol\theta_{ij}(\boldsymbol\sigma)+\boldsymbol\theta_{ij}(\boldsymbol\eta)-\boldsymbol\tau\Big) -\psi(\boldsymbol\theta_{ij}(\boldsymbol\sigma))-\psi(\boldsymbol\theta_{ij}(\boldsymbol\eta))+\psi(\boldsymbol\tau)\\  \nonumber
&=&2\sum_{|\alpha|=2}\frac{\partial^{\alpha}\psi}{\alpha!}(\boldsymbol\tau)\left(\frac{\boldsymbol\tau_{\boldsymbol d}}{\sqrt{n}}\right)^{\alpha}+14\sum_{|\alpha|=4}\frac{\partial^{\alpha}\psi}{\alpha!}(\boldsymbol\tau)\left(\frac{\boldsymbol\tau_{\boldsymbol d}}{\sqrt{n}}\right)^{\alpha}-6\sum_{|\alpha|=3}\frac{\partial^{\alpha}\psi}{\alpha!}(\boldsymbol\tau)\left(\frac{\boldsymbol\tau_{\boldsymbol d}}{\sqrt{n}}\right)^{\alpha}+O\left(\frac{1}{n^2\sqrt{n}}\right).\\ \label{exp6}
\end{eqnarray}

For $\sigma_i\sigma_j=1$ and $\eta_i\eta_j=-1$ or $\sigma_i\sigma_j=-1$ and $\eta_i\eta_j=1$, one has
\begin{eqnarray}\nonumber
&&\psi\Big(\boldsymbol\theta_{ij}(\boldsymbol\sigma)+\boldsymbol\theta_{ij}(\boldsymbol\eta)-\boldsymbol\tau\Big) -\psi(\boldsymbol\theta_{ij}(\boldsymbol\sigma))-\psi(\boldsymbol\theta_{ij}(\boldsymbol\eta))+\psi(\boldsymbol\tau)\\ \nonumber
&=&\psi\left(\tau_1,\tau_2\right)-\psi\left(\tau_1+\frac{\tau_1d_1}{\sqrt{n}},\tau_2+\frac{\tau_2d_2}{\sqrt{n}}\right)-\psi\left(\tau_1-\frac{\tau_1d_1}{\sqrt{n}},\tau_2-\frac{\tau_2d_2}{\sqrt{n}}\right)+\psi(\tau_1,\tau_2)\\ \label{exp7}
&=&-2\sum_{|\alpha|=2}\frac{\partial^{\alpha}\psi}{\alpha!}(\boldsymbol\tau)\left(\frac{\boldsymbol\tau_{\boldsymbol d}}{\sqrt{n}}\right)^{\alpha}-2\sum_{|\alpha|=4}\frac{\partial^{\alpha}\psi}{\alpha!}(\boldsymbol\tau)\left(\frac{\boldsymbol\tau_{\boldsymbol d}}{\sqrt{n}}\right)^{\alpha}+O\left(\frac{1}{n^2\sqrt{n}}\right).
\end{eqnarray}

For $\sigma_i\sigma_j=-1$ and $\eta_i\eta_j=-1$, the following equations are true.
\begin{eqnarray}\nonumber
&&\psi\Big(\boldsymbol\theta_{ij}(\boldsymbol\sigma)+\boldsymbol\theta_{ij}(\boldsymbol\eta)-\boldsymbol\tau\Big) -\psi(\boldsymbol\theta_{ij}(\boldsymbol\sigma))-\psi(\boldsymbol\theta_{ij}(\boldsymbol\eta))+\psi(\boldsymbol\tau)\\ \nonumber
&=&\psi\left(\tau_1+\frac{2\tau_1d_1}{\sqrt{n}},\tau_2+\frac{2\tau_2d_2}{\sqrt{n}}\right)-2\psi\left(\tau_1+\frac{\tau_1d_1}{\sqrt{n}},\tau_2+\frac{\tau_2d_2}{\sqrt{n}}\right)+\psi(\tau_1,\tau_2)\\ \nonumber
&=&2\sum_{|\alpha|=2}\frac{\partial^{\alpha}\psi}{\alpha!}(\boldsymbol\tau)\left(\frac{\boldsymbol\tau_{\boldsymbol d}}{\sqrt{n}}\right)^{\alpha}+14\sum_{|\alpha|=4}\frac{\partial^{\alpha}\psi}{\alpha!}(\boldsymbol\tau)\left(\frac{\boldsymbol\tau_{\boldsymbol d}}{\sqrt{n}}\right)^{\alpha}+6\sum_{|\alpha|=3}\frac{\partial^{\alpha}\psi}{\alpha!}(\boldsymbol\tau)\left(\frac{\boldsymbol\tau_{\boldsymbol d}}{\sqrt{n}}\right)^{\alpha}+O\left(\frac{1}{n^2\sqrt{n}}\right).\\ \label{exp8}
\end{eqnarray}

Let
\[
s_+=\#\{(i,j): i<j, \sigma_i\sigma_j\eta_i\eta_j=+1\},\ \ \ \
s_-=\#\{(i,j): i<j, \sigma_i\sigma_j\eta_i\eta_j=-1\},
\]
and $\rho=\frac{1}{n}\sum_{i=1}^n\sigma_i\eta_i$. Then
$s_+=\frac{n^2}{4}(1+\rho^2)-\frac{n}{2}$
and $s_-=\frac{n^2}{4}(1-\rho^2)$. By (\ref{exp2}) and (\ref{exp6})-(\ref{exp8}), we get
\begin{eqnarray}\nonumber
&&\mathbb{E}[L_n^2]=(1+o(1))\mathbb{E}_{\boldsymbol\sigma,\boldsymbol\eta}\Bigg[\exp\Bigg\{-s_-\left(2\sum_{|\alpha|=2}\frac{\partial^{\alpha}\psi}{\alpha!}(\boldsymbol\tau)\left(\frac{\boldsymbol\tau_{\boldsymbol d}}{\sqrt{n}}\right)^{\alpha}+2\sum_{|\alpha|=4}\frac{\partial^{\alpha}\psi}{\alpha!}(\boldsymbol\tau)\left(\frac{\boldsymbol\tau_{\boldsymbol d}}{\sqrt{n}}\right)^{\alpha}\right)\Bigg\}\\  \nonumber
&&\times\exp\left\{s_+\left(2\sum_{|\alpha|=2}\frac{\partial^{\alpha}\psi}{\alpha!}(\boldsymbol\tau)\left(\frac{\boldsymbol\tau_{\boldsymbol d}}{\sqrt{n}}\right)^{\alpha}+14\sum_{|\alpha|=4}\frac{\partial^{\alpha}\psi}{\alpha!}(\boldsymbol\tau)\left(\frac{\boldsymbol\tau_{\boldsymbol d}}{\sqrt{n}}\right)^{\alpha}\right)\right\}\\ \nonumber 
&&\times\exp\left\{6\sum_{|\alpha|=3}\frac{\partial^{\alpha}\psi}{\alpha!}(\boldsymbol\tau)\left(\frac{\boldsymbol\tau_{\boldsymbol d}}{\sqrt{n}}\right)^{\alpha}\sum_{i<j}\left(\frac{(1-\sigma_i\sigma_j)}{2}\frac{(1-\eta_i\eta_j)}{2}-\frac{(1+\sigma_i\sigma_j)}{2}\frac{(1+\eta_i\eta_j)}{2}\right)\right\}
\Bigg]\\ \nonumber
&=&(1+o(1))\mathbb{E}_{\boldsymbol\sigma,\boldsymbol\eta}\Bigg[\exp\Bigg\{2\frac{s_+-s_-}{n}\sum_{|\alpha|=2}\frac{\partial^{\alpha}\psi}{\alpha!}(\boldsymbol\tau)\left(\boldsymbol\tau_{\boldsymbol d}\right)^{\alpha}\Bigg\}\exp\Bigg\{\frac{14s_+-2s_-}{n^2}\sum_{|\alpha|=4}\frac{\partial^{\alpha}\psi}{\alpha!}(\boldsymbol\tau)\left(\boldsymbol\tau_{\boldsymbol d}\right)^{\alpha}\Bigg\}\\ \nonumber 
&&\times\exp\left\{-3\sum_{|\alpha|=3}\frac{\partial^{\alpha}\psi}{\alpha!}(\boldsymbol\tau)\left(\boldsymbol\tau_{\boldsymbol d}\right)^{\alpha}\left[\frac{1}{n\sqrt{n}}\left(\sum_{i=1}^n\sigma_i\right)^2+\frac{1}{n\sqrt{n}}\left(\sum_{i=1}^n\eta_i\right)^2\right]\right\}\Bigg]\exp\Bigg\{-\sum_{|\alpha|=2}\frac{\partial^{\alpha}\psi}{\alpha!}(\boldsymbol\tau)\left(\boldsymbol\tau_{\boldsymbol d}\right)^{\alpha}\Bigg\}\\ \label{exp9} 
&=&(1+o(1))\mathbb{E}_{\boldsymbol\sigma,\boldsymbol\eta}\big[Y_n\big]\exp\Bigg\{-\sum_{|\alpha|=2}\frac{\partial^{\alpha}\psi}{\alpha!}(\boldsymbol\tau)\left(\boldsymbol\tau_{\boldsymbol d}\right)^{\alpha}+3\sum_{|\alpha|=4}\frac{\partial^{\alpha}\psi}{\alpha!}(\boldsymbol\tau)\left(\boldsymbol\tau_{\boldsymbol d}\right)^{\alpha}\Bigg\}, 
\end{eqnarray}
where
\begin{eqnarray}
Y_n&=&\exp\Bigg\{n\rho^2\sum_{|\alpha|=2}\frac{\partial^{\alpha}\psi}{\alpha!}(\boldsymbol\tau)\left(\boldsymbol\tau_{\boldsymbol d}\right)^{\alpha}\Bigg\}\exp\Bigg\{4\rho^2\sum_{|\alpha|=4}\frac{\partial^{\alpha}\psi}{\alpha!}(\boldsymbol\tau)\left(\boldsymbol\tau_{\boldsymbol d}\right)^{\alpha}\Bigg\}\\ \nonumber
&&\times\exp\left\{-3\sum_{|\alpha|=3}\frac{\partial^{\alpha}\psi}{\alpha!}(\boldsymbol\tau)\left(\boldsymbol\tau_{\boldsymbol d}\right)^{\alpha}\left[\frac{1}{n\sqrt{n}}\left(\sum_{i=1}^n\sigma_i\right)^2+\frac{1}{n\sqrt{n}}\left(\sum_{i=1}^n\eta_i\right)^2\right]\right\}.
\end{eqnarray}

Next we find the limit of $\mathbb{E}_{\boldsymbol\sigma,\boldsymbol\eta}\big[Y_n\big]$.
Note that $ \exp\Bigg\{n\rho^2\sum_{|\alpha|=2}\frac{\partial^{\alpha}\psi}{\alpha!}(\boldsymbol\tau)\left(\boldsymbol\tau_{\boldsymbol d}\right)^{\alpha}\Bigg\}$ is uniformly integrable if $2\sum_{|\alpha|=2}\frac{\partial^{\alpha}\psi}{\alpha!}(\boldsymbol\tau)\left(\boldsymbol\tau_{\boldsymbol d}\right)^{\alpha}=\boldsymbol\tau_{\boldsymbol d}^TD^2\psi(\boldsymbol\tau)\boldsymbol\tau_{\boldsymbol d}<1$.
 Besides, $n\rho^2$ converges in law to chi-square distribution with degree of freedom one.  Hence 
 \begin{equation}\label{exp10}
 \mathbb{E}\left(\exp\Bigg\{\frac{n\rho^2}{2}\boldsymbol\tau_{\boldsymbol d}^TD^2\psi(\boldsymbol\tau)\boldsymbol\tau_{\boldsymbol d}\Bigg\}\right)=(1+o(1))\left(1-\boldsymbol\tau_{\boldsymbol d}^TD^2\psi(\boldsymbol\tau)\boldsymbol\tau_{\boldsymbol d}\right)^{-\frac{1}{2}}.
 \end{equation}

Let $M=32\sum_{|\alpha|=3}\frac{|\partial^{\alpha}\psi|}{\alpha!}(\boldsymbol\tau)\left(\boldsymbol\tau_{\boldsymbol d}\right)^{\alpha}$. Decompose $Y_n$ as follows.
\[Y_n=Y_{n,1}+Y_{n,2}+Y_{n,3}+Y_{n,4},\]
where
\begin{eqnarray*}
Y_{n,1}&=&Y_nI\Big[\frac{(\sum_{i}\sigma_i)^2}{n\sqrt{n}}\leq M,\frac{(\sum_{i}\eta_i)^2}{n\sqrt{n}}\leq M\Big],\\
Y_{n,2}&=&Y_nI\Big[\frac{(\sum_{i}\sigma_i)^2}{n\sqrt{n}}> M,\frac{(\sum_{i}\eta_i)^2}{n\sqrt{n}}\leq M\Big],
\\
Y_{n,3}&=&Y_nI\Big[\frac{(\sum_{i}\sigma_i)^2}{n\sqrt{n}}> M,\frac{(\sum_{i}\eta_i)^2}{n\sqrt{n}}> M\Big],\\
Y_{n,4}&=&Y_nI\Big[\frac{(\sum_{i}\sigma_i)^2}{n\sqrt{n}}\leq M,\frac{(\sum_{i}\eta_i)^2}{n\sqrt{n}}> M\Big].
\end{eqnarray*}

Let $p=1+\frac{1-\boldsymbol\tau_{\boldsymbol d}^TD^2\psi(\boldsymbol\tau)\boldsymbol\tau_{\boldsymbol d}}{2\boldsymbol\tau_{\boldsymbol d}^TD^2\psi(\boldsymbol\tau)\boldsymbol\tau_{\boldsymbol d} }$. Then $p>1$ and $ p\boldsymbol\tau_{\boldsymbol d}^TD^2\psi(\boldsymbol\tau)\boldsymbol\tau_{\boldsymbol d}<1$. In this case, by (\ref{exp10}),
\[\mathbb{E}_{\boldsymbol\sigma,\boldsymbol\eta}\Big[Y_{n1}^p\Big]=O(1)\mathbb{E}\left(\exp\Bigg\{\frac{n\rho^2}{2}p\boldsymbol\tau_{\boldsymbol d}^TD^2\psi(\boldsymbol\tau)\boldsymbol\tau_{\boldsymbol d}\Bigg\}\right)=O(1),\]
which implies $Y_{n1}$ is uniformly intregrable. 
Note that $Y_{n1}$ converges in distribution to $\exp\Big\{\frac{\chi^2_1}{2}\boldsymbol\tau_{\boldsymbol d}^TD^2\psi(\boldsymbol\tau)\boldsymbol\tau_{\boldsymbol d}\Big\}$. Hence
 \begin{equation}\label{ep10}
 \mathbb{E}\left(Y_{n1}\right)=(1+o(1))\left(1-\boldsymbol\tau_{\boldsymbol d}^TD^2\psi(\boldsymbol\tau)\boldsymbol\tau_{\boldsymbol d}\right)^{-\frac{1}{2}}.
 \end{equation}

Next we show $\mathbb{E}_{\boldsymbol\sigma,\boldsymbol\eta}\Big[Y_{nt}\Big]=o(1)$ for $t=2,3,4$.
By Bernstein inequality, for  $n>M^2$, we have
\begin{equation}\label{exp11}
\mathbb{P}\left(\frac{(\sum_{i}\sigma_i)^2}{n\sqrt{n}}> M\right)=\mathbb{P}\left(\Big|\frac{1}{n}\sum_{i}\sigma_i\Big|>\sqrt{\frac{M}{\sqrt{n}}}\right)\leq 2\exp\left(-\frac{\sqrt{n}M}{4}\right).
\end{equation}
Since $\sigma_i\eta_i$ and $\sigma_i$ are independent, then $\rho$ and $\sum_{i=1}^n\sigma_i$ are independent. By the fact that $\frac{(\sum_{i}\sigma_i)^2}{n\sqrt{n}}\leq \sqrt{n}$, $\frac{(\sum_{i}\eta_i)^2}{n\sqrt{n}}\leq \sqrt{n}$, (\ref{exp10})  and (\ref{exp11}), one has
\begin{eqnarray}\nonumber
&&\mathbb{E}_{\boldsymbol\sigma,\boldsymbol\eta}\Big[Y_{n2}\Big]\\ \nonumber
&\leq&O(1)\exp\left\{6\sqrt{n}\sum_{|\alpha|=3}\frac{|\partial^{\alpha}\psi|}{\alpha!}(\boldsymbol\tau)\left(\boldsymbol\tau_{\boldsymbol d}\right)^{\alpha}\right\}\mathbb{E}_{\boldsymbol\sigma,\boldsymbol\eta}\Bigg[\exp\Bigg\{n\rho^2\sum_{|\alpha|=2}\frac{\partial^{\alpha}\psi}{\alpha!}(\boldsymbol\tau)\left(\boldsymbol\tau_{\boldsymbol d}\right)^{\alpha}\Bigg\}I\Big[\frac{(\sum_{i}\sigma_i)^2}{n\sqrt{n}}> M\Big]\Bigg]\\ \nonumber
&=&O(1)\exp\left\{6\sqrt{n}\sum_{|\alpha|=3}\frac{|\partial^{\alpha}\psi|}{\alpha!}(\boldsymbol\tau)\left(\boldsymbol\tau_{\boldsymbol d}\right)^{\alpha}\right\}\mathbb{E}_{\boldsymbol\sigma,\boldsymbol\eta}\Bigg[\exp\Bigg\{n\rho^2\sum_{|\alpha|=2}\frac{\partial^{\alpha}\psi}{\alpha!}(\boldsymbol\tau)\left(\boldsymbol\tau_{\boldsymbol d}\right)^{\alpha}\Bigg\}\Bigg]\mathbb{P}\left(\frac{(\sum_{i}\sigma_i)^2}{n\sqrt{n}}> M\right)\\ \nonumber
&\leq&2O(1)\exp\left\{6\sqrt{n}\sum_{|\alpha|=3}\frac{|\partial^{\alpha}\psi|}{\alpha!}(\boldsymbol\tau)\left(\boldsymbol\tau_{\boldsymbol d}\right)^{\alpha}-\frac{\sqrt{n}M}{4}\right\}\mathbb{E}_{\boldsymbol\sigma,\boldsymbol\eta}\Bigg[\exp\Bigg\{n\rho^2\sum_{|\alpha|=2}\frac{\partial^{\alpha}\psi}{\alpha!}(\boldsymbol\tau)\left(\boldsymbol\tau_{\boldsymbol d}\right)^{\alpha}\Bigg\}\Bigg]\\ \label{exp12}
&=&o(1).
\end{eqnarray}
Similarly, $\mathbb{E}_{\boldsymbol\sigma,\boldsymbol\eta}\Big[Y_{nt}\Big]=o(1)$ for $t=3,4$. Then by (\ref{exp9}) and (\ref{exp12}), 
\begin{equation}\label{sedm}
\mathbb{E}(L_n^2)=(1+o(1))\exp\Bigg\{-\sum_{|\alpha|=2}\frac{\partial^{\alpha}\psi}{\alpha!}(\boldsymbol\tau)\left(\boldsymbol\tau_{\boldsymbol d}\right)^{\alpha}+3\sum_{|\alpha|=4}\frac{\partial^{\alpha}\psi}{\alpha!}(\boldsymbol\tau)\left(\boldsymbol\tau_{\boldsymbol d}\right)^{\alpha}\Bigg\}\left(1-\boldsymbol\tau_{\boldsymbol d}^TD^2\psi(\boldsymbol\tau)\boldsymbol\tau_{\boldsymbol d}\right)^{-\frac{1}{2}}.
\end{equation}
Hence, $\mathbb{E}(L_n^2)=O(1)$ under $H_0$, which implies any test is inconsistent.

{\bf (II)}. The proof strategy is to show the WSLMC test has asymptotic power one. 
To this end, we prove that $\mathcal{Z}_n=O_P(1)$ under $H_0$ and   $\mathcal{Z}_n=\frac{(-1)^k}{\sqrt{2k}}\left[\boldsymbol\tau_{\boldsymbol d}^TD^2\psi(\boldsymbol\tau)\boldsymbol\tau_{\boldsymbol d}\right]^{\frac{k}{2}}(1+o(1))+o_P(1)$ under $H_1$.

By the property of exponential family, the mean and covariance of $\boldsymbol T(X)$ are equal to
\[\mathbb{E}(\boldsymbol T(X))=D\psi,\ \ \ \ \ \  Cov(\boldsymbol T(X))=D^2\psi.\]
Under $H_0$, $A_{ij} (1\leq i<j\leq n)$ are independent. Hence, $\mathbb{E}(\mathcal{Z}_n)=0$. Note that
\begin{eqnarray}\nonumber
&&\sum_{1\leq i_1<i_2<\dots<i_k\leq n}\mathbb{E}\left[\prod_{t=1}^k\frac{\boldsymbol\tau_{\boldsymbol d}^T\left(\boldsymbol T(A_{i_ti_{t+1}})-D\psi(\boldsymbol\tau)\right)}{\sqrt{n}}\right]^2\\ \nonumber
&=&\sum_{1\leq i_1<i_2<\dots<i_k\leq n}\prod_{t=1}^k\frac{\boldsymbol\tau_{\boldsymbol d}^T \mathbb{E}\left[(\boldsymbol T(A_{i_ti_{t+1}})-D\psi(\boldsymbol\tau))(\boldsymbol T(A_{i_ti_{t+1}})-D\psi(\boldsymbol\tau))^T\right]\boldsymbol\tau_{\boldsymbol d}}{n}\\ \label{cycl1}
&=&\binom{n}{k}\frac{\left[\boldsymbol\tau_{\boldsymbol d}^TD^2\psi(\boldsymbol\tau)\boldsymbol\tau_{\boldsymbol d}\right]^{k}}{n^{k}}.
\end{eqnarray}
For $(i_1,i_2,\dots,i_k)\in \mathcal{I}_n$ and $(j_1,j_2,\dots,j_k)\in \mathcal{I}_n$, if $(i_1,i_2,\dots,i_k)\neq(j_1,j_2,\dots,j_k)$, then \[\mathbb{E}\left[\prod_{t=1}^k\boldsymbol\tau_{\boldsymbol d}^T\left(\boldsymbol T(A_{i_ti_{t+1}})-D\psi(\boldsymbol\tau)\right)\prod_{t=1}^k\boldsymbol\tau_{\boldsymbol d}^T\left(\boldsymbol T(A_{j_tj_{t+1}})-D\psi(\boldsymbol\tau)\right)\right]=0.\]
For  $k=\log\log\log n$, it is easy to verify $\binom{n}{k}=\frac{n^{k}}{k!}(1+o(1))$.
Hence, by (\ref{cycl1}), we have
\[\mathbb{E}\mathcal{Z}_n^2=\frac{\frac{(k-1)!}{2}\binom{n}{k}\frac{\left[\boldsymbol\tau_{\boldsymbol d}^TD^2\psi(\boldsymbol\tau)\boldsymbol\tau_{\boldsymbol d}\right]^{k}}{n^k}}{\frac{1}{2k}\left[\boldsymbol\tau_{\boldsymbol d}^TD^2\psi(\boldsymbol\tau)\boldsymbol\tau_{\boldsymbol d}\right]^k}=1+o(1),\]
which implies $\mathcal{Z}_n=O_p(1)$ under $H_0$.

Next, we show $\mathcal{Z}_n=\frac{1}{\sqrt{2k}}\left[\boldsymbol\tau_{\boldsymbol d}^TD^2\psi(\boldsymbol\tau)\boldsymbol\tau_{\boldsymbol d}\right]^{\frac{k}{2}}+o_p(1)$ under $H_1$. Given $\boldsymbol \sigma$, the mean and covariance of $\boldsymbol T(A_{ij})$ are equal to
\[\mathbb{E}(\boldsymbol T(A_{ij}))=D\psi(\boldsymbol\theta_{ij}(\boldsymbol\sigma)),\hskip 1cm Cov(\boldsymbol T(A_{ij}))=D^2\psi(\boldsymbol\theta_{ij}(\boldsymbol\sigma)),\]
where $\boldsymbol\theta_{ij}(\boldsymbol\sigma)$ is defined in (\ref{theta1}). Let $\lambda_i\in\{0,1\}, (1\leq i\leq k)$ and $\boldsymbol\lambda=(\lambda_1,\dots,\lambda_k)$. Denote $|\boldsymbol\lambda|=\lambda_1+\lambda_2+\dots+\lambda_k$. Then 
\begin{eqnarray}\nonumber
&&\mathcal{Z}_n\\ \nonumber
&=&\frac{\sum_{(i_1,i_2,\dots,i_k)\in\mathcal{I}_n}\prod_{t=1}^k\boldsymbol\tau_{\boldsymbol d}^T\left(\boldsymbol T(A_{i_ti_{t+1}})-D\psi(\boldsymbol\theta_{i_ti_{t+1}}(\boldsymbol\sigma))+D\psi(\boldsymbol\theta_{i_ti_{t+1}}(\boldsymbol\sigma))-D\psi(\boldsymbol\tau)\right)}{\sqrt{\frac{1}{2k}n^k\left[\boldsymbol\tau_{\boldsymbol d}^TD^2\psi(\boldsymbol\tau)\boldsymbol\tau_{\boldsymbol d}\right]^k}}\\ \nonumber
&=&\frac{\sum_{(i_1,i_2,\dots,i_k)\in\mathcal{I}_n}\prod_{t=1}^k\boldsymbol\tau_{\boldsymbol d}^T\left(D\psi(\boldsymbol\theta_{i_ti_{t+1}}(\boldsymbol\sigma))-D\psi(\boldsymbol\tau)\right)}{\sqrt{\frac{1}{2k}n^k\left[\boldsymbol\tau_{\boldsymbol d}^TD^2\psi(\boldsymbol\tau)\boldsymbol\tau_{\boldsymbol d}\right]^k}}\\  \nonumber
&+&\frac{\sum_{(i_1,i_2,\dots,i_k)\in\mathcal{I}_n}\sum_{|\boldsymbol\lambda|>0}\prod_{t=1}^k\Big[\boldsymbol\tau_{\boldsymbol d}^T\left(\boldsymbol T(A_{i_ti_{t+1}})-D\psi(\boldsymbol\theta_{i_ti_{t+1}}(\boldsymbol\sigma))\Big)\Big]^{\lambda_t}\Big[\boldsymbol\tau_{\boldsymbol d}^T\Big(D\psi(\boldsymbol\theta_{i_ti_{t+1}}(\boldsymbol\sigma))-D\psi(\boldsymbol\tau)\right)\Big]^{1-\lambda_t}}{\sqrt{\frac{1}{2k}n^k\left[\boldsymbol\tau_{\boldsymbol d}^TD^2\psi(\boldsymbol\tau)\boldsymbol\tau_{\boldsymbol d}\right]^k}}\\ \label{cycl3}
&=&R_1+R_2. 
\end{eqnarray}
Next we show the first term $R_1$ in (\ref{cycl3}) is the leading term.

By Taylor expansion, we get
\begin{eqnarray*}\nonumber
\frac{\partial\psi}{\partial \theta_1}(\boldsymbol\theta_{ij}(\boldsymbol\sigma))
&=&\frac{\partial\psi}{\partial \theta_1}(\boldsymbol\tau)+\frac{\partial^2\psi}{\partial\theta_1^2}(\boldsymbol\tau)\frac{(-1)\tau_1d_1}{\sqrt{n}}\sigma_i\sigma_j+\frac{\partial^2\psi}{\partial\theta_1\partial\theta_2}(\boldsymbol\tau)\frac{(-1)\tau_2d_2}{\sqrt{n}}\sigma_i\sigma_j+O\left(\frac{1}{n}\right),\\ \label{cycl2}
\frac{\partial\psi}{\partial \theta_2}(\boldsymbol\theta_{ij}(\boldsymbol\sigma))&=&\frac{\partial\psi}{\partial \theta_2}(\boldsymbol\tau)+\frac{\partial^2\psi}{\partial\theta_2\partial\theta_1}(\boldsymbol\tau)\frac{(-1)\tau_1d_1}{\sqrt{n}}\sigma_i\sigma_j+\frac{\partial^2\psi}{\partial\theta_2^2}(\boldsymbol\tau)\frac{(-1)\tau_2d_2}{\sqrt{n}}\sigma_i\sigma_j+O\left(\frac{1}{n}\right).
\end{eqnarray*}
Hence $\boldsymbol\tau_{\boldsymbol d}^T\Big(D\psi(\boldsymbol\theta_{ij}(\boldsymbol\sigma))-D\psi(\boldsymbol\tau)\Big)=-\frac{\boldsymbol\tau_{\boldsymbol d}^TD^2\psi(\boldsymbol\tau)\boldsymbol\tau_{\boldsymbol d}}{\sqrt{n}}\sigma_i\sigma_j+O\left(\frac{1}{n}\right)$ and
\begin{equation}\label{cycl4}
\prod_{t=1}^k\boldsymbol\tau_{\boldsymbol d}^T\left(D\psi(\boldsymbol\theta_{i_ti_{t+1}}(\boldsymbol\sigma))-D\psi(\boldsymbol\tau)\right)=\prod_{t=1}^k\left(-\frac{\boldsymbol\tau_{\boldsymbol d}^TD^2\psi(\boldsymbol\tau)\boldsymbol\tau_{\boldsymbol d}}{\sqrt{n}}+O\left(\frac{1}{n}\right)\right)\sigma_{i_t}\sigma_{i_{t+1}}=\frac{[-\boldsymbol\tau_{\boldsymbol d}^TD^2\psi(\boldsymbol\tau)\boldsymbol\tau_{\boldsymbol d}]^k}{\sqrt{n^k}}(1+o(1)).
\end{equation}
Then (\ref{cycl4}) implies $R_1=(1+o(1))(-1)^k\sqrt{\frac{1}{2k}[\boldsymbol\tau_{\boldsymbol d}^TD^2\psi(\boldsymbol\tau)\boldsymbol\tau_{\boldsymbol d}]^k}$.

Next, we prove $R_2=o_p(1)$. Since $\boldsymbol\tau_{\boldsymbol d}$ is a vector of constants, then for a large constant $C>0$,
\[\mathbb{E}\left[\boldsymbol\tau_{\boldsymbol d}\Big(\boldsymbol T(A_{i_ti_{t+1}})-D\psi(\boldsymbol\theta_{i_ti_{t+1}}(\boldsymbol\sigma))\Big)\right]^2\leq C,\]
\[\left|\boldsymbol\tau_{\boldsymbol d}\Big(D\psi(\boldsymbol\theta_{ij}(\boldsymbol\sigma))-D\psi(\boldsymbol\tau)\Big)\right|\leq\frac{C}{\sqrt{n}}.\]
Recall that $A_{ij}(1\leq i<j\leq n)$ are independent conditional on $\boldsymbol\sigma$.
Then fixing a $\boldsymbol\lambda$ with $|\boldsymbol\lambda|=l>0$, we have
\begin{eqnarray*}
&&\mathbb{E}\Bigg[\frac{\sum_{i_1<i_2<\dots<i_k}\prod_{t=1}^k\Big[\boldsymbol\tau_{\boldsymbol d}\left(\boldsymbol T(A_{i_ti_{t+1}})-D\psi(\boldsymbol\theta_{i_ti_{t+1}}(\boldsymbol\sigma))\Big)\Big]^{\lambda_t}\Big[\boldsymbol\tau_{\boldsymbol d}\Big(D\psi(\boldsymbol\theta_{i_ti_{t+1}}(\boldsymbol\sigma))-D\psi(\boldsymbol\tau)\right)\Big]^{1-\lambda_t}}{\sqrt{\frac{1}{2k}n^k\left[\boldsymbol\tau_{\boldsymbol d}D^2\psi(\boldsymbol\tau)\boldsymbol\tau_{\boldsymbol d}\right]^k}}\Bigg]^2\\
&=&\frac{1}{\frac{1}{2k}n^k\left[\boldsymbol\tau_{\boldsymbol d}D^2\psi(\boldsymbol\tau)\boldsymbol\tau_{\boldsymbol d}\right]^k}\sum_{\substack{i_1<i_2<\dots<i_k\\ j_1<j_2<\dots<j_k\\
(j_t,j_{t+1})=(i_t,i_{t+1})\ if\ \lambda_t=1}}\prod_{t=1}^k\mathbb{E}\Bigg\{\Big[\boldsymbol\tau_{\boldsymbol d}\Big(\boldsymbol T(A_{i_ti_{t+1}})-D\psi(\boldsymbol\theta_{i_ti_{t+1}}(\boldsymbol\sigma))\Big)\Big]^{2\lambda_t}\\
&&\times \Big[\boldsymbol\tau_{\boldsymbol d}\Big(D\psi(\boldsymbol\theta_{i_ti_{t+1}}(\boldsymbol\sigma))-D\psi(\boldsymbol\tau)\Big)\Big]^{1-\lambda_t}\Big[\boldsymbol\tau_{\boldsymbol d}\Big(D\psi(\boldsymbol\theta_{j_tj_{t+1}}(\boldsymbol\sigma))-D\psi(\boldsymbol\tau)\Big)\Big]^{1-\lambda_t}\Bigg\}\\
&\leq&\frac{2k}{n^k\left[\boldsymbol\tau_{\boldsymbol d}D^2\psi(\boldsymbol\tau)\boldsymbol\tau_{\boldsymbol d}\right]^k}n^{2k-s}C^{l}\frac{C^{2k-2l}}{\sqrt{n}^{2k-2l}}\\
&=&\frac{2kC^{2k-l}}{n^{s-l}\left[\boldsymbol\tau_{\boldsymbol d}D^2\psi(\boldsymbol\tau)\boldsymbol\tau_{\boldsymbol d}\right]^k}\leq \frac{2kC^{2k}}{n\left[\boldsymbol\tau_{\boldsymbol d}D^2\psi(\boldsymbol\tau)\boldsymbol\tau_{\boldsymbol d}\right]^k}.
\end{eqnarray*}
Here, $s$ is the number of distinct nodes that any $l$ edges on the cycle $i_1,i_2,\dots,i_k$ have and hence $s\geq l+1$. 
Note that there are $2^k-1$ possible choices of $\boldsymbol\lambda$ such that $|\boldsymbol\lambda|=l>0$. Since $k!\leq k^{k}$ and $k=\log\log\log n$, then
\[\mathbb{E}(R_2^2)\leq \frac{((k-1)!)^2}{4}2^k\frac{2kC^{2k}}{n\left[\boldsymbol\tau_{\boldsymbol d}D^2\psi(\boldsymbol\tau)\boldsymbol\tau_{\boldsymbol d}\right]^k}=o(1).\]
Then the proof is complete.

\qed

\subsection{Proof of Theorem \ref{expoloss}}

The proof strategy is similar to that of Theorem \ref{theorem:1}. The likelihood ratio $L_n$ is equal to
\begin{eqnarray*}
L_n=\frac{1}{2^n}\sum_{\boldsymbol\sigma\in\{\pm\}^n}\left[\prod_{1\leq i<j\leq n}\frac{p_{ij}(\sigma)^{\tilde{A}_{ij}}\left(1-p_{ij}(\sigma)\right)^{1-\tilde{A}_{ij}}}{p_0^{\tilde{A}_{ij}}(1-p_0)^{1-\tilde{A}_{ij}}}\right].
\end{eqnarray*}
Let $\boldsymbol\eta$ be an independent copy of $\boldsymbol\sigma$. Then the second moment under $H_0$ is
\begin{eqnarray}\nonumber
\mathbb{E}(L_n^2)&=&\mathbb{E}_{\boldsymbol\sigma,\boldsymbol\eta}\left[\prod_{1\leq i<j\leq n}\left(\frac{p_{ij}(\boldsymbol\sigma)p_{ij}(\boldsymbol\eta)}{p_0}+\frac{(1-p_{ij}(\boldsymbol\sigma))(1-p_{ij}(\boldsymbol\eta))}{1-p_0}\right)\right]\\ \label{elnh0}
&=&\mathbb{E}_{\boldsymbol\sigma,\boldsymbol\eta}\left[\prod_{1\leq i<j\leq n}\frac{p_0-p_0[p_{ij}(\boldsymbol\sigma)+p_{ij}(\boldsymbol\eta)]+p_{ij}(\boldsymbol\sigma)p_{ij}(\boldsymbol\eta)}{p_0(1-p_0)}\right].
\end{eqnarray}

By Taylor expansion, we have
\begin{eqnarray}\nonumber
\psi\left(\boldsymbol\tau-\frac{\boldsymbol\tau_{\boldsymbol d}} {\sqrt{n}}\sigma_i\sigma_j\right)&=&\psi(\boldsymbol\tau)-\frac{\boldsymbol\tau_{\boldsymbol d}}{\sqrt{n}}D\psi(\boldsymbol\tau)\sigma_i\sigma_j+\sum_{|\alpha|=2}\frac{\partial^{\alpha}\psi}{\alpha!}(\boldsymbol\tau)\left(\frac{\boldsymbol\tau_{\boldsymbol d}}{\sqrt{n}}\right)^{\alpha}+O\left(\frac{1}{n\sqrt{n}}\right)\sigma_i\sigma_j.
\end{eqnarray}
Hence,
\begin{eqnarray}\nonumber
&&p_{ij}(\boldsymbol\sigma)-p_0\\ \nonumber
&=&\int_{t_0}^{\infty}h(x)e^{\boldsymbol\theta(\sigma_i,\sigma_j)^T\boldsymbol T(x)-\psi(\boldsymbol\theta(\sigma_i,\sigma_j))}-h(x)e^{\boldsymbol\tau^T\boldsymbol T(x)-\psi(\boldsymbol\tau)}dx\\ \nonumber
&=&\int_{t_0}^{\infty}h(x)e^{\boldsymbol\tau^T\boldsymbol T(x)-\psi(\boldsymbol\tau)}\Bigg[\exp\Bigg\{\frac{\boldsymbol\tau_{\boldsymbol d}}{\sqrt{n}}\big(-\boldsymbol T(x)+D\psi(\boldsymbol\tau)\big)\sigma_i\sigma_j\\ \nonumber
&&-\sum_{|\alpha|=2}\frac{\partial^{\alpha}\psi}{\alpha!}(\boldsymbol\tau)\left(\frac{\boldsymbol\tau_{\boldsymbol d}}{\sqrt{n}}\right)^{\alpha}+O\left(\frac{1}{n\sqrt{n}}\right)\sigma_i\sigma_j\Bigg\}-1\Bigg]dx\\ \nonumber
&=&\int_{t_0}^{\infty}h(x)e^{\boldsymbol\tau^T\boldsymbol T(x)-\psi(\boldsymbol\tau)}\Bigg[\Bigg(\frac{\boldsymbol\tau_{\boldsymbol d}}{\sqrt{n}}\big(-\boldsymbol T(x)+D\psi(\boldsymbol\tau)\big)\sigma_i\sigma_j-\sum_{|\alpha|=2}\frac{\partial^{\alpha}\psi}{\alpha!}(\boldsymbol\tau)\left(\frac{\boldsymbol\tau_{\boldsymbol d}}{\sqrt{n}}\right)^{\alpha}+O\left(\frac{1}{n\sqrt{n}}\right)\sigma_i\sigma_j\Bigg)\\ \nonumber
&&+\Bigg(\frac{\boldsymbol\tau_{\boldsymbol d}^T}{\sqrt{n}}\big(-\boldsymbol T(x)+D\psi(\boldsymbol\tau)\big)\sigma_i\sigma_j-\sum_{|\alpha|=2}\frac{\partial^{\alpha}\psi}{\alpha!}(\boldsymbol\tau)\left(\frac{\boldsymbol\tau_{\boldsymbol d}}{\sqrt{n}}\right)^{\alpha}+O\left(\frac{1}{n\sqrt{n}}\right)\sigma_i\sigma_j\Bigg)^2\\ \nonumber
&&+\Bigg(\frac{\boldsymbol\tau_{\boldsymbol d}^T}{\sqrt{n}}\big(-\boldsymbol T(x)+D\psi(\boldsymbol\tau)\big)\sigma_i\sigma_j-\sum_{|\alpha|=2}\frac{\partial^{\alpha}\psi}{\alpha!}(\boldsymbol\tau)\left(\frac{\boldsymbol\tau_{\boldsymbol d}}{\sqrt{n}}\right)^{\alpha}+O\left(\frac{1}{n\sqrt{n}}\right)\sigma_i\sigma_j\Bigg)^3\Bigg]dx+O\left(\frac{1}{n^2}\right)\\ \nonumber
&=&\frac{\sigma_i\sigma_j\boldsymbol\tau_{\boldsymbol d}^T\boldsymbol a(t_0)}{\sqrt{n}}+\frac{b(t_0)}{n}+O\left(\frac{1}{n\sqrt{n}}\right)\sigma_i\sigma_j,
\end{eqnarray}
where
\[b(t_0)=\int_{t_0}^{\infty}h(x)e^{\boldsymbol\tau^T\boldsymbol T(x)-\psi(\boldsymbol\tau)}\left(\Big(\boldsymbol\tau_{\boldsymbol d}^T\big(-\boldsymbol T(x)+D\psi(\boldsymbol\tau)\Big)^2-\sum_{|\alpha|=2}\frac{\partial^{\alpha}\psi}{\alpha!}(\boldsymbol\tau)\left(\boldsymbol\tau_{\boldsymbol d}\right)^{\alpha}\right).\]
Then we get
\[p_{ij}(\boldsymbol\sigma)=p_0+\frac{\sigma_i\sigma_j\boldsymbol\tau_{\boldsymbol d}^T\boldsymbol a(t_0)}{\sqrt{n}}+\frac{b(t_0)}{n}+O\left(\frac{1}{n\sqrt{n}}\right)\sigma_i\sigma_j+O\left(\frac{1}{n^2}\right),\]
and then
\begin{eqnarray}\nonumber
&&p_0-p_0[p_{ij}(\boldsymbol\sigma)+p_{ij}(\boldsymbol\eta)]+p_{ij}(\boldsymbol\sigma)p_{ij}(\boldsymbol\eta)\\ \nonumber
&=&p_0-p_0\Bigg(2p_0+\frac{2b(t_0)}{n}+\frac{(\sigma_i\sigma_j+\eta_i\eta_j)\boldsymbol\tau_{\boldsymbol d}^T\boldsymbol a(t_0)}{\sqrt{n}}+O\left(\frac{1}{n\sqrt{n}}\right)(\sigma_i\sigma_j+\eta_i\eta_j)+O\left(\frac{1}{n^2}\right)\Bigg)\\ \nonumber
&&+p_0^2+\frac{2p_0b(t_0)}{n}+\frac{(\sigma_i\sigma_j+\eta_i\eta_j)p_0\boldsymbol\tau_{\boldsymbol d}^T\boldsymbol a(t_0)}{\sqrt{n}}+\frac{(\boldsymbol\tau_{\boldsymbol d}^T\boldsymbol a(t_0))^2}{n}\sigma_i\sigma_j\eta_i\eta_j+O\left(\frac{1}{n\sqrt{n}}\right)(\sigma_i\sigma_j+\eta_i\eta_j)\\ \label{p0pij}
&=&p_0(1-p_0)+\frac{(\boldsymbol\tau_{\boldsymbol d}^T\boldsymbol a(t_0))^2}{n}\sigma_i\sigma_j\eta_i\eta_j+O\left(\frac{1}{n\sqrt{n}}\right)(\sigma_i\sigma_j+\eta_i\eta_j).
\end{eqnarray}

Plugging (\ref{p0pij}) into (\ref{elnh0}) yields
\begin{eqnarray*}\nonumber
\mathbb{E}(L_n^2)
&=&\mathbb{E}_{\boldsymbol\sigma,\boldsymbol\eta}\left[\prod_{1\leq i<j\leq n}\frac{p_0-p_0[p_{ij}(\boldsymbol\sigma)+p_{ij}(\boldsymbol\eta)]+p_{ij}(\boldsymbol\sigma)p_{ij}(\boldsymbol\eta)}{p_0(1-p_0)}\right]\\ \nonumber
&=&\mathbb{E}_{\boldsymbol\sigma,\boldsymbol\eta}\left[\prod_{1\leq i<j\leq n}\left(1+\frac{(\boldsymbol\tau_{\boldsymbol d}^T\boldsymbol a(t_0))^2}{np_0(1-p_0)}\sigma_i\sigma_j\eta_i\eta_j+O\left(\frac{1}{n\sqrt{n}}\right)(\sigma_i\sigma_j+\eta_i\eta_j)\right)\right]\\ \nonumber
&=&O(1)\mathbb{E}_{\boldsymbol\sigma,\boldsymbol\eta}\Bigg[\exp\Bigg(\frac{(\boldsymbol\tau_{\boldsymbol d}^T\boldsymbol a(t_0))^2}{p_0(1-p_0)}\frac{\sum_{i<j}\sigma_i\sigma_j\eta_i\eta_j}{n}+O\Big(1\Big)\frac{\sum_{i<j}(\sigma_i\sigma_j+\eta_i\eta_j)}{n\sqrt{n}}\Bigg)\Bigg]\\
&=&O(1)\mathbb{E}_{\boldsymbol\sigma,\boldsymbol\eta}\Bigg[\exp\Bigg(\frac{(\boldsymbol\tau_{\boldsymbol d}^T\boldsymbol a(t_0))^2}{2p_0(1-p_0)}\frac{\big(\sum_{i}\sigma_i\eta_i\big)^2-n}{n}+O\Big(1\Big)\frac{\big(\sum_{i}\sigma_i\big)^2+\big(\sum_{i}\eta_i\big)^2-2n}{2n\sqrt{n}}\Bigg)\Bigg]
\end{eqnarray*}
If $\frac{(\boldsymbol\tau_{\boldsymbol d}^T\boldsymbol a(t_0))^2}{p_0(1-p_0)}<1$, then $\exp\Bigg(\frac{(\boldsymbol\tau_{\boldsymbol d}^T\boldsymbol a(t_0))^2}{2p_0(1-p_0)}\frac{\big(\sum_{i}\sigma_i\eta_i\big)^2}{n}\Bigg)$ is uniformly integrable and 
\[\mathbb{E}_{\boldsymbol\sigma,\boldsymbol\eta}\Bigg[\exp\exp\Bigg(\frac{(\boldsymbol\tau_{\boldsymbol d}^T\boldsymbol a(t_0))^2}{2p_0(1-p_0)}\frac{\big(\sum_{i}\sigma_i\eta_i\big)^2}{n}\Bigg)\Bigg]=(1+o(1))\left(1-\frac{(\boldsymbol\tau_{\boldsymbol d}^T\boldsymbol a(t_0))^2}{p_0(1-p_0)}\right)^{-\frac{1}{2}}.\]
By a similar truncation technique in the proof of Theorem \ref{theorem:1}, we conclude $\mathbb{E}(L_n^2)=O(1)$ if $\frac{(\boldsymbol\tau_{\boldsymbol d}^T\boldsymbol a(t_0))^2}{p_0(1-p_0)}<1$.

{\bf (II)}. By a similar proof of (II) of Theorem \ref{theorem:1},
 $\mathcal{Z}_n=O_P(1)$ under $H_0$ and   $\mathcal{Z}_n=\frac{(-1)^k}{\sqrt{2k}}\left[\frac{(\boldsymbol\tau_{\boldsymbol d}^T\boldsymbol a(t_0))^2}{p_0(1-p_0)}\right]^{\frac{k}{2}}(1+o(1))+o_P(1)$ under $H_1$. Then the proof is complete.

\qed

\subsection{Proof of Theorem \ref{contiguity}}

Firstly we recall the contiguity theorem (\cite{B18}) and several useful lemmas.

\begin{Proposition}\label{conti}
Let $\mathbb{P}_n$  and $\mathbb{Q}_n$ be two sequences of probability measures and $X_{ni}$ be random variables on the same sample space. Then $\mathbb{P}_n$  and $\mathbb{Q}_n$ are mutually contiguous if the following conditions hold.\\
$i)$.  $\mathbb{P}_n<<\mathbb{Q}_n$ and $\mathbb{Q}_n<<\mathbb{P}_n$.\\
$ii)$. For any fixed $m\geq3$, $X_{ni} (3\leq i\leq m)$ jointly converges in distribution to $Z_i$ with $ Z_i\sim N(0,2i)$ under $\mathbb{P}_n$ and $Y_i$ with $ Y_i\sim N(t^{\frac{i}{2}},2i)$ ($|t|<1$) under  $\mathbb{Q}_n$ respectively.\\
$iii)$. $Z_i$ and $Y_i$ $(3\leq i\leq m)$ are independent.\\
$iv)$. 
\[\mathbb{E}_{\mathbb{P}_n}\left[\frac{d\mathbb{Q}_n}{d\mathbb{P}_n}\right]^2=(1+o(1))\exp\left(-\frac{t}{2}-\frac{t^2}{4}\right)\frac{1}{\sqrt{1-t}}.\]
\end{Proposition}

The following two lemmas are well-known.

\begin{Lemma}\label{cntlem1}
Let $Y_{n,1},\dots,Y_{n,m}$ be $m$ random variables. Then $Y_{n,1},\dots,Y_{n,m}$ jointly converges in distribution to $Z_1,\dots,Z_m$ if the following conditions hold.\\
i). For any fixed $k$ and  $\lambda_1+\dots+\lambda_m=k$ with integers $\lambda_t\geq0$, $\mathbb{E}[Y_{n,1}^{\lambda_1}\dots,Y_{n,m}^{\lambda_m}]=(1+o(1))\mathbb{E}[Z_{1}^{\lambda_1}\dots,Z_{m}^{\lambda_m}]$.\\
ii). 
\[\sum_{t=1}^{\infty}\left(\lim_{n\rightarrow\infty}\mathbb{E}(Y_{n,l}^{2t})\right)^{-\frac{1}{2t}},\ \ 1\leq l\leq m.\]
\end{Lemma}

\begin{Lemma}\label{cntlem2}
Let $Y_1,\dots,Y_m$ follow a $m$-variate distribution $F$ with mean 0 and covariance $\Sigma$. Then $F$ is Gaussian distribution if and only if for even $l$,
\[\mathbb{E}[X_1X_2\dots X_l]=\sum_{\eta}\prod_{i=1}^{\frac{l}{2}}\mathbb{E}[X_{\eta(i,1)}X_{\eta(i,2)}]\]
and $\mathbb{E}[X_1X_2\dots X_l]=0$ for odd $l$. Here $X_i\in\{Y_1,Y_2,\dots, Y_m\}$, $\eta$ is a partition of $\{1,2,\dots,l\}$ into $\frac{l}{2}$ equal-size subsets and $\eta(i,t)$ is the $t$th element of $i$th subset.
\end{Lemma}

Given integer $k\geq3$, let $\mathcal{J}_{n,k}=\{(i_1,i_2,\dots,i_k)|i_1,\dots,i_k: distinct\}$ and define
\begin{eqnarray*}
\mathcal{U}_{n,k}=\frac{\sum_{(i_1,i_2,\dots,i_k)\in\mathcal{J}_{n,k}}\prod_{t=1}^k\boldsymbol\tau_{\boldsymbol d}^T\left(\boldsymbol T(A_{i_ti_{t+1}})-D\psi(\boldsymbol\tau)\right)}{\sqrt{n^k\left[\boldsymbol\tau_{\boldsymbol d}^TD^2\psi(\boldsymbol\tau)\boldsymbol\tau_{\boldsymbol d}\right]^k}}.
\end{eqnarray*}

\begin{Proposition}\label{jtpd}
The following results hold.\\
(a). For fixed integers $3\leq k_1<k_2<\dots<k_m$, $\frac{\mathcal{U}_{n,k_t}}{\sqrt{2k_t}}$ ($1\leq t\leq m$) converges jointly in distribution to the standard $m$-variate normal distribution under $H_0$.\\
(b). For fixed integers $3\leq k_1<k_2<\dots<k_m$, $\frac{\mathcal{U}_{n,k_t}-\left[\boldsymbol\tau_{\boldsymbol d}^TD^2\psi(\boldsymbol\tau)\boldsymbol\tau_{\boldsymbol d}\right]^{k_t}}{\sqrt{2k_t}}$ ($1\leq t\leq m$) converges jointly in distribution to the standard $m$-variate normal distribution under $H_1$.
\end{Proposition}

\noindent
\textit{Proof of Proposition \ref{jtpd}}. By the proof of Theorem \ref{theorem:1}, it is easy to get that $Var(\mathcal{U}_{n,k})=2k(1+o(1))$ for each fixed integer $k\geq3$ under $H_0$ or $H_1$. By Lemma \ref{cntlem1} and Lemma \ref{cntlem2}, to prove (a) or (b), it suffices to prove that
\begin{equation}\label{cnt1}
\mathbb{E}[\mathcal{U}_{n,k_1}\mathcal{U}_{n,k_2}]=o(1),\ \ k_1\neq k_2,
\end{equation}
and
\begin{equation}\label{cnt2}
\mathbb{E}[\mathcal{V}_{n,1}\dots\mathcal{V}_{n,l}]
=
\begin{cases}
\sum_{\eta}\prod_{i=1}^{\frac{l}{2}}\mathbb{E}[Z_{\eta(i,1)}Z_{\eta(i,2)}],\ \text{ even $l$},\\      
 0,  \text{ odd $l$},
\end{cases}
\end{equation}
where $\mathcal{V}_{n,i}\in\{\mathcal{U}_{n,k_1},\dots,\mathcal{U}_{n,k_m}\}$.

(a). Under $H_0$, $A_{ij}$ $(1\leq i<j\leq n)$ are independent. Hence it is easy to verify that (\ref{cnt1}) holds. Next we prove (\ref{cnt2}). Let $\lambda_i\in\{k_1,k_2,\dots,k_m\}$ be the length of cycles in $\mathcal{V}_{n,i}$. Let $\boldsymbol i_t=(i_t^{(1)},\dots,i_t^{(\lambda_t)})\in\mathcal{J}_{n,\lambda_t}$ and $V_{\boldsymbol i_t}=\prod_{j=1}^{\lambda_t}\boldsymbol\tau_{\boldsymbol d}^T\left(\boldsymbol T(A_{i_t^{(j)}i_t^{(j+1)}})-D\psi(\boldsymbol\tau)\right)$. Then
\begin{equation}\label{cnt3}
\mathbb{E}[\mathcal{V}_{n,1}\dots\mathcal{V}_{n,l}]=\frac{\sum_{\boldsymbol i_1,\dots,\boldsymbol i_l}\mathbb{E}(V_{\boldsymbol i_1}\dots V_{\boldsymbol i_l})}{\left(n\left[\boldsymbol\tau_{\boldsymbol d}^TD^2\psi(\boldsymbol\tau)\boldsymbol\tau_{\boldsymbol d}\right]\right)^{\frac{\sum_{i=1}^l\lambda_i}{2}}}.
\end{equation}
Each edge in $V_{\boldsymbol i_t}$ must be traversed at least twice, otherwise  $\mathbb{E}(V_{\boldsymbol i_1}\dots V_{\boldsymbol i_l})=0$. Hence any node pair $(i_{t_1}^{j_1},i_{t_1}^{j_1+1})\in \boldsymbol i_{t_1}$ must be equal to at least one other pair $(i_{t_2}^{j_2},i_{t_2}^{j_2+1})\in \boldsymbol i_{t_2}$ with $t_1\neq t_2$. Then $\boldsymbol i_{t}$ ($1\leq t\leq l$) are partitioned into $s$ disjoint groups with each group containing at least 2 elements. Clearly $s\leq\frac{l}{2}$. If $s<\frac{l}{2}$, then $\boldsymbol i_{t}$ ($1\leq t\leq l$) has at most $\frac{\sum_i\lambda_i}{2}-1$ distinct nodes. Hence by (\ref{cnt3}), we have
\[\mathbb{E}[\mathcal{V}_{n,1}\dots\mathcal{V}_{n,l}]=O\left(\frac{n^{\frac{\sum_i\lambda_i}{2}-1}}{n^{\frac{\sum_{i=1}^l\lambda_i}{2}}}\right)=o(1).\]
Note that for odd $l$, $s<\frac{l}{2}$ always holds. When $l$ is even and $s=\frac{l}{2}$, it is easy to check that (\ref{cnt2}) holds.

(b). By the proof of Part (II) of Theorem \ref{theorem:1} and part (a) above, the proof is straightforward. Hence we omit it.

\qed

\textit{Proof of Theorem \ref{contiguity}}: 
By (\ref{sedm}), if (\ref{ODE}) holds, then $iv)$ of Proposition \ref{conti} holds with $t=\boldsymbol\tau_{\boldsymbol d}D^2\psi(\boldsymbol\tau)\boldsymbol\tau_{\boldsymbol d}$. The proof is straightforward based on Proposition \ref{conti} and Proposition \ref{jtpd}.

\qed

\subsection{Proof of Theorem \ref{prac1} and Theorem \ref{prac2}}

\noindent
{\it Proof of Theorem \ref{prac1}:} Suppose $H_0$ holds. Let $\Sigma=Cov(\boldsymbol M(A_{12}))$.  Then for $k=\log\log\log n$,
$(\boldsymbol 1^T\boldsymbol S^21)^k$ converges to $(\boldsymbol 1^T\Sigma1)^k$ in probability. Let
\begin{eqnarray*}
\mathcal{Y}_n=\frac{\sum_{(i_1,i_2,\dots,i_k)\in\mathcal{I}_n}\prod_{t=1}^k\boldsymbol 1^T\left(\boldsymbol M(A_{i_ti_{t+1}})-\overline{\boldsymbol M}( \boldsymbol A)\right)}{\sqrt{\frac{(k-1)!}{2}\binom{n}{k}\big(\boldsymbol 1^T\boldsymbol \Sigma1\big)^k}}.
\end{eqnarray*}
We only need to prove $\mathcal{Y}_n$ converges to the standard normal distribution. Let $\lambda_t\in\{0,1\}$, $\boldsymbol\lambda=(\lambda_1,\dots,\lambda_k)$ and $|\boldsymbol\lambda|=\lambda_1+\lambda_2+\dots+\lambda_k$. Note that
\begin{eqnarray}\nonumber
&&\sum_{i_1<i_2<\dots<i_k}\prod_{t=1}^k\boldsymbol 1^T\left(\boldsymbol M(A_{i_ti_{t+1}})-\overline{\boldsymbol M}( \boldsymbol A)\right)\\ \nonumber
&=&\sum_{i_1<i_2<\dots<i_k}\prod_{t=1}^k\boldsymbol 1^T\left(\boldsymbol M(A_{i_ti_{t+1}})-\boldsymbol\mu+\boldsymbol\mu-\overline{\boldsymbol M}( \boldsymbol A)\right)\\ \nonumber
&=&\sum_{i_1<i_2<\dots<i_k}\prod_{t=1}^k\boldsymbol 1^T\left(\boldsymbol M(A_{i_ti_{t+1}})-\boldsymbol\mu\right)+\sum_{i_1<i_2<\dots<i_k}\sum_{|\boldsymbol\lambda|<k}\prod_{t=1}^k[\boldsymbol 1^T\left(\boldsymbol M(A_{i_ti_{t+1}})-\boldsymbol\mu\right)]^{\lambda_t}[\boldsymbol 1^T(\boldsymbol\mu-\overline{\boldsymbol M}( \boldsymbol A))]^{1-\lambda_t}\\ \label{power1}
&=&R_1+R_2.
\end{eqnarray}
Next we show $\mathbb{E}(R_1^2)=\binom{n}{k}\big(\boldsymbol 1^T\boldsymbol \Sigma1\big)^k$ and $R_2=o_p\left(\sqrt{\binom{n}{k}\big(\boldsymbol 1^T\boldsymbol \Sigma1\big)^k}\right)$. 

Under $H_0$, $A_{ij} (1\leq i<j\leq n)$ are independent and $\mathbb{E}(\boldsymbol M(A_{ij}))=\boldsymbol\mu$. Hence
\begin{eqnarray}\label{pow2}
\mathbb{E}\Bigg[\sum_{i_1<i_2<\dots<i_k}\prod_{t=1}^k\boldsymbol 1^T\left(\boldsymbol M(A_{i_ti_{t+1}})-\boldsymbol\mu\right)\Bigg]^2=\sum_{i_1<i_2<\dots<i_k}\prod_{t=1}^k\mathbb{E}\Big[\boldsymbol 1^T\left(\boldsymbol M(A_{i_ti_{t+1}})-\boldsymbol\mu\right)\Big]^2=\binom{n}{k}\big(\boldsymbol 1^T\boldsymbol \Sigma1\big)^k.
\end{eqnarray}

For a given $\boldsymbol\lambda$ with $|\boldsymbol\lambda|<k$, we have
\begin{eqnarray*}
&&\sum_{i_1<i_2<\dots<i_k}\prod_{t=1}^k[\boldsymbol 1^T\left(\boldsymbol M(A_{i_ti_{t+1}})-\boldsymbol\mu\right)]^{\lambda_t}[\boldsymbol 1^T(\boldsymbol\mu-\overline{\boldsymbol M}( \boldsymbol A))]^{1-\lambda_t}\\
&=&\prod_{t=1}^k[\boldsymbol 1^T(\boldsymbol\mu-\overline{\boldsymbol M}( \boldsymbol A))]^{1-\lambda_t}\sum_{i_1<i_2<\dots<i_k}\prod_{t=1}^k[\boldsymbol 1^T\left(\boldsymbol M(A_{i_ti_{t+1}})-\boldsymbol\mu\right)]^{\lambda_t}.
\end{eqnarray*}
Note that
\begin{eqnarray*}
&&\mathbb{E}\Bigg[\sum_{i_1<i_2<\dots<i_k}\prod_{t=1}^k[\boldsymbol 1^T\left(\boldsymbol M(A_{i_ti_{t+1}})-\boldsymbol\mu\right)]^{\lambda_t}\Bigg]^2\\
&=&\sum_{\substack{i_1<i_2<\dots<i_k\\ j_1<j_2<\dots<j_k\\(i_t,i_{t+1})=(j_t,j_{t+1})\ if\ \lambda_t=1 }}\prod_{t=1}^k\mathbb{E}[\boldsymbol 1^T\left(\boldsymbol M(A_{i_ti_{t+1}})-\boldsymbol\mu\right)]^{2\lambda_t}\\
&\leq&n^{2k-s}\big(\boldsymbol 1^T\boldsymbol \Sigma1\big)^{|\boldsymbol\lambda|},
\end{eqnarray*}
where $k>s\geq|\boldsymbol\lambda|+1$. Since $\boldsymbol\mu-\overline{\boldsymbol M}( \boldsymbol A)=O_p\left(\frac{1}{n}\right)$, then
\begin{equation}\label{pow3}
R_2=O_p\left(2^kn^{|\boldsymbol\lambda|-\frac{s}{2}}\big(\boldsymbol 1^T\boldsymbol \Sigma1\big)^{\frac{|\boldsymbol\lambda|}{2}}\right)=o_p\left(\sqrt{\binom{n}{k}\big(\boldsymbol 1^T\boldsymbol \Sigma1\big)^k}\right).
\end{equation}

By (\ref{power1}), (\ref{pow2}) and (\ref{pow3}), we conclude that
\begin{eqnarray*}
\mathcal{Y}_n=\mathcal{X}_n+o_P(1),
\end{eqnarray*}
where
\[\mathcal{X}_n=\frac{\sum_{(i_1,i_2,\dots,i_k)\in\mathcal{I}_n}\prod_{t=1}^k\boldsymbol 1^T\left(\boldsymbol M(A_{i_ti_{t+1}})-\boldsymbol\mu\right)}{\sqrt{\frac{(k-1)!}{2}\binom{n}{k}\big( \boldsymbol1^T\boldsymbol \Sigma \boldsymbol1\big)^k}}.\]

Next, we use the method of moment to prove $\mathcal{X}_n$ converges in distribution to the standard normal distribution. To this end, we will show $\mathbb{E}(\mathcal{X}_n^r)=o(1)$ for odd $r$ and $\mathbb{E}(\mathcal{X}_n^r)=(r-1)!!+o(1)$ for even $r$.

Clearly, $\mathbb{E}(\mathcal{X}_n)=0$. The second moment of $\mathcal{X}_n$ is
\begin{eqnarray*}
\mathbb{E}(\mathcal{X}_n^2)&=&\frac{\sum_{\substack{(i_1,i_2,\dots,i_k)\in\mathcal{I}_n\\(j_1,j_2,\dots,j_k)\in\mathcal{I}_n}}\mathbb{E}\Big[\prod_{t=1}^k\boldsymbol 1^T\left(\boldsymbol M(A_{i_ti_{t+1}})-\boldsymbol\mu\right)\boldsymbol 1^T\left(\boldsymbol M(A_{j_tj_{t+1}})-\boldsymbol\mu\right)\Big]}{\frac{(k-1)!}{2}\binom{n}{k}\big(\boldsymbol 1^T\boldsymbol \Sigma1\big)^k}\\
&=&\frac{\sum_{\substack{(i_1,i_2,\dots,i_k)\in\mathcal{I}_n}}\prod_{t=1}^k\mathbb{E}\Big[\boldsymbol 1^T\left(\boldsymbol M(A_{i_ti_{t+1}})-\boldsymbol\mu\right)\Big]^2}{\frac{(k-1)!}{2}\binom{n}{k}\big(\boldsymbol 1^T\boldsymbol \Sigma1\big)^k}=1.
\end{eqnarray*}

Fix a positive integer $r\geq3$.
For convenience, let $\boldsymbol i_k^{(t)}=(i_1^{(t)},i_2^{(t)},\dots,i_k^{(t)})$ denote a circular permutation of $k$ distinct nodes for each $t\in\{1,2,\dots, r\}$. Then
 
\begin{eqnarray*}
\mathbb{E}(\mathcal{X}_n^r)&=&\frac{\sum_{\boldsymbol i_k^{(1)},\dots,\boldsymbol i_k^{(r)}\in\mathcal{I}_n}\mathbb{E}\Big[\prod_{v=1}^r\prod_{t=1}^k\boldsymbol 1^T\left(\boldsymbol M(A_{i_t^{(v)}i_{t+1}^{(v)}})-\boldsymbol\mu\right)\Big]}{\left[\frac{(k-1)!}{2}\binom{n}{k}\big(\boldsymbol 1^T\boldsymbol \Sigma1\big)^k\right]^{\frac{r}{2}}}.
\end{eqnarray*}

If there are two indexes $i_t^{(v)},i_{t+1}^{(v)}$ such that $(i_t^{(v)},i_{t+1}^{(v)})$ is different from any other pairs, then \[\mathbb{E}\Big[\prod_{v=1}^r\prod_{t=1}^k\boldsymbol 1^T\left(\boldsymbol M(A_{i_t^{(v)}i_{t+1}^{(v)}})-\boldsymbol\mu\right)\Big]=0.\]
Hence, any $\boldsymbol i_k^{(v_1)}\in \mathcal{I}_n$ has to be equal to at least one $\boldsymbol i_k^{(v_2)}\in \mathcal{I}_n$ for $v_2\neq v_1$. Then there exist $s (s\leq\frac{r}{2})$ integers $\lambda_i, (\lambda_i\geq2,i=1,2,\dots,s)$ such that $\lambda_1+\lambda_2+\dots+\lambda_s=r$ and
\begin{eqnarray*}
\mathbb{E}(\mathcal{X}_n^r)&=&C_r\frac{\sum_{\boldsymbol i_k^{(1)},\dots,\boldsymbol i_k^{(s)}\in\mathcal{I}_n}\mathbb{E}\Big[\prod_{v=1}^s\prod_{t=1}^k[\boldsymbol 1^T\left(\boldsymbol M(A_{i_t^{(v)}i_{t+1}^{(v)}})-\boldsymbol\mu\right)]^{\lambda_v}\Big]}{\Big[\frac{(k-1)!}{2}\binom{n}{k}\big(\boldsymbol 1^T\boldsymbol \Sigma1\big)^k\Big]^{\frac{r}{2}}},
\end{eqnarray*}
where $C_r$ is a constant dependent on $r$. Note that all the moments of $A_{12}$ are finite. By repeatedly using Cauchy-Schwarz inequality, we have
\[\mathbb{E}\Big[\prod_{v=1}^s\prod_{t=1}^k\big[\boldsymbol 1^T\left(\boldsymbol M(A_{i_t^{(v)}i_{t+1}^{(v)}})-\boldsymbol\mu\right)\big]^{\lambda_v}\Big]\leq C^{kr},\]
for a large constant $C$. If there exists $\lambda_v\geq3$, that is, $s<\frac{r}{2}$, then
\begin{eqnarray*}
\mathbb{E}(\mathcal{X}_n^r)&\leq&C_r\frac{n^{ks}C^{kr}}{\Big[\frac{(k-1)!}{2}\binom{n}{k}\big(\boldsymbol 1^T\boldsymbol \Sigma1\big)^k\Big]^{\frac{r}{2}}}=\frac{C_r}{\big(\boldsymbol 1^T\boldsymbol \Sigma1\big)^{\frac{kr}{2}}}\frac{(2k)^{\frac{r}{2}}C^{kr}}{n^{k(\frac{r}{2}-s)}}=o(1),
\end{eqnarray*}
noting that $k=\log\log\log n$. If $r$ is an odd number, $s<\frac{r}{2}$ holds and hence $\mathbb{E}(\mathcal{X}_n^r)=o(1)$.

Next we assume $r$ is even and $s=\frac{r}{2}$. Note that there are  $(r-1)!!$ ways to partition $r$ distinct numbers into $\frac{r}{2}$ pairs. Then $C_r=(r-1)!!$ and 
\begin{eqnarray*}
\mathbb{E}(\mathcal{X}_n^r)&=&(r-1)!!\frac{\sum_{\boldsymbol i_k^{(1)},\dots,\boldsymbol i_k^{(s)}\in\mathcal{I}_n}\mathbb{E}\Big[\prod_{v=1}^s\prod_{t=1}^k[\boldsymbol 1^T\left(\boldsymbol M(A_{i_t^{(v)}i_{t+1}^{(v)}})-\boldsymbol\mu\right)]^{2}\Big]}{\Big[\frac{(k-1)!}{2}\binom{n}{k}\big(\boldsymbol 1^T\boldsymbol \Sigma1\big)^k\Big]^{\frac{r}{2}}}.
\end{eqnarray*}
If there are two $\boldsymbol i_k^{(v_1)}$ and $\boldsymbol i_k^{(v_2)}$ ($v_1\neq v_2$) have at least a common vertex, then $\mathbb{E}(\mathcal{X}_n^r)=o(1)$. Hence 
\begin{eqnarray*}
\mathbb{E}(\mathcal{X}_n^r)&=&(r-1)!!\frac{\sum_{\substack{\boldsymbol i_k^{(1)},\dots,\boldsymbol i_k^{(s)}\in\mathcal{I}_n\\ \boldsymbol i_k^{(1)},\dots,\boldsymbol i_k^{(s)}:disjoint }}\prod_{v=1}^s\mathbb{E}\Big[\prod_{t=1}^k[\boldsymbol 1^T\left(\boldsymbol M(A_{i_t^{(v)}i_{t+1}^{(v)}})-\boldsymbol\mu\right)]^{2}\Big]}{\Big[\frac{(k-1)!}{2}\binom{n}{k}\big(\boldsymbol 1^T\boldsymbol \Sigma1\big)^k\Big]^{\frac{r}{2}}}=(r-1)!!.
\end{eqnarray*}
Then the proof is complete.

\qed

\noindent
{\it Proof of Theorem \ref{prac2}:} We find the order of each term in (\ref{power1}). Given $\boldsymbol\sigma$, the mean of $\boldsymbol M(A_{ij})$ is equal to $\boldsymbol\mu_{ij}=\boldsymbol\mu+\frac{\boldsymbol\mu_{\boldsymbol d}}{\sqrt{n}}\sigma_i\sigma_j$. Hence,
\begin{eqnarray}\nonumber
&&\sum_{i_1<i_2<\dots<i_k}\prod_{t=1}^k\boldsymbol 1^T\left(\boldsymbol M(A_{i_ti_{t+1}})-\boldsymbol\mu\right)\\ \nonumber
&=&\sum_{i_1<i_2<\dots<i_k}\prod_{t=1}^k\boldsymbol 1^T\left(\boldsymbol M(A_{i_ti_{t+1}})-\boldsymbol\mu_{i_ti_{t+1}}+\boldsymbol\mu_{i_ti_{t+1}}-\boldsymbol\mu\right)\\ \nonumber
&=&\sum_{i_1<i_2<\dots<i_k}\prod_{t=1}^k\boldsymbol 1^T(\boldsymbol\mu_{i_ti_{t+1}}-\boldsymbol\mu)
+\sum_{i_1<i_2<\dots<i_k}\prod_{t=1}^k\boldsymbol 1^T\left(\boldsymbol M(A_{i_ti_{t+1}})-\boldsymbol\mu_{i_ti_{t+1}}\right)
\\  \nonumber
&&+
\sum_{i_1<i_2<\dots<i_k}\sum_{0<|\boldsymbol\lambda|<k}\prod_{t=1}^k[\boldsymbol 1^T\left(\boldsymbol M(A_{i_ti_{t+1}})-\boldsymbol\mu_{i_ti_{t+1}}\right)]^{\lambda_t}[\boldsymbol 1^T(\boldsymbol\mu_{i_ti_{t+1}}-\boldsymbol\mu)]^{1-\lambda_t}\\ \label{power3}
&=&\binom{n}{k}\left(\frac{\boldsymbol 1^T\boldsymbol\mu_{\boldsymbol d}}{\sqrt{n}}\right)^k+(a)+(b).
\end{eqnarray}
Note that the conditional second moment $\boldsymbol\Sigma_{ij}$ of $\boldsymbol M(A_{ij})$ exists and is a function of $\boldsymbol\mu_{ij}=\boldsymbol\mu+\frac{\boldsymbol\mu_{\boldsymbol d}}{\sqrt{n}}\sigma_i\sigma_j$. By Taylor expansion, it follows that $\boldsymbol 1^T\boldsymbol\Sigma_{ij}\boldsymbol 1=\boldsymbol 1^T\boldsymbol\Sigma\boldsymbol 1+O\left(\frac{\max_{t}|d_t|}{\sqrt{n}}\right)$. Hence,
\begin{eqnarray}\nonumber
\mathbb{E}((a)^2)&=&\sum_{i_1<i_2<\dots<i_k}\mathbb{E}\left(\prod_{t=1}^k\mathbb{E}\left[\left(\boldsymbol 1^T\left(\boldsymbol M(A_{i_ti_{t+1}})-\boldsymbol\mu_{i_ti_{t+1}}\right)\right)^2\Big|\boldsymbol\sigma\right]\right)\\ \nonumber
&=&\sum_{i_1<i_2<\dots<i_k}\mathbb{E}\left[ \prod_{t=1}^k\left(\boldsymbol 1^T\boldsymbol\Sigma\boldsymbol 1+O\left(\frac{\max_{t}|d_t|}{\sqrt{n}}\right)\right)\right]\\ \label{power2}
&=&\binom{n}{k}\left(\boldsymbol 1^T\boldsymbol\Sigma\boldsymbol 1\right)^k(1+o(1)).
\end{eqnarray}
For $s$ with $1+|\boldsymbol\lambda|\leq s\leq 2|\boldsymbol\lambda|$, the order of $(b)$ is bounded by
\begin{equation}\label{power4}
2^k\left(\frac{\boldsymbol 1^T\boldsymbol\mu_{\boldsymbol d}}{\sqrt{n}}\right)^{k-|\boldsymbol\lambda|}\sqrt{n^{2k-s}\left(\boldsymbol 1^T\boldsymbol\Sigma\boldsymbol 1+O\left(\frac{\max_{t}|d_t|}{\sqrt{n}}\right)\right)^{|\boldsymbol\lambda|}}=o\left(\binom{n}{k}\left(\frac{\boldsymbol 1^T\boldsymbol\mu_{\boldsymbol d}}{\sqrt{n}}\right)^k\right).
\end{equation}

Hence by (\ref{power3}), (\ref{power2}) and (\ref{power4}), we have
\[\frac{\sum_{(i_1,i_2,\dots,i_k)\in\mathcal{I}_n}\prod_{t=1}^k\boldsymbol 1^T\left(\boldsymbol M(A_{i_ti_{t+1}})-\boldsymbol\mu\right)}{\sqrt{\frac{(k-1)!}{2}\binom{n}{k}\big(\boldsymbol1^T\boldsymbol \Sigma\boldsymbol 1\big)^k}}=\frac{1}{\sqrt{2k}}\left(\frac{\boldsymbol 1^T\boldsymbol\mu_{\boldsymbol d}}{\sqrt{\boldsymbol1^T\boldsymbol \Sigma\boldsymbol 1}}\right)^k(2+o_p(1)).\]

Given $\boldsymbol\sigma$, we have
\begin{eqnarray}\nonumber
\boldsymbol 1^T\left(\boldsymbol\mu-\overline{\boldsymbol M}( \boldsymbol A)\right)&=&\frac{1}{\binom{n}{2}}\sum_{1\leq i<j\leq n}\boldsymbol 1^T\left(\boldsymbol\mu-\boldsymbol\mu_{ij}+\boldsymbol\mu_{ij}-\boldsymbol M( \boldsymbol A_{ij})\right)\\  \nonumber
&=&\frac{1}{\binom{n}{2}}\sum_{1\leq i<j\leq n} \frac{-\boldsymbol1^T\boldsymbol\mu_{\boldsymbol d}}{\sqrt{n}}\sigma_i\sigma_j+\frac{1}{\binom{n}{2}}\sum_{1\leq i<j\leq n}\boldsymbol 1^T\left(\boldsymbol\mu_{ij}-\boldsymbol M( \boldsymbol A_{ij})\right)\\  \nonumber
&=&\frac{-\boldsymbol1^T\boldsymbol\mu_{\boldsymbol d}}{\sqrt{n}}\frac{\left(\sum_{i=1}^n\sigma_i\right)^2-n}{n(n-1)}+O_p\left(\sqrt{\frac{\sum_{1\leq i<j\leq n}\boldsymbol1^T\boldsymbol\Sigma_{ij}\boldsymbol1}{\binom{n}{2}^2}}\right)\\
&=&O_p\left(\frac{\boldsymbol1^T\boldsymbol\mu_{\boldsymbol d}}{\sqrt{n}n}+\frac{1}{n}\right)=O_p\left(\frac{1}{n}\right),
\end{eqnarray}
and
\[\boldsymbol1^T\left(\boldsymbol M(A_{i_ti_{t+1}})-\boldsymbol\mu\right)=\boldsymbol1^T\left(\boldsymbol M(A_{i_ti_{t+1}})-\boldsymbol\mu_{ij}\right)+\frac{\boldsymbol1^T\boldsymbol\mu_{\boldsymbol d}}{\sqrt{n}}\sigma_i\sigma_j=\boldsymbol1^T\left(\boldsymbol M(A_{i_ti_{t+1}})-\boldsymbol\mu_{ij}\right)+o_p(1).\]
For a given $\boldsymbol\lambda$, 
\begin{eqnarray}\nonumber
&&\mathbb{E}\left[\sum_{i_1<i_2<\dots<i_k}\prod_{t=1}^k[\boldsymbol 1^T\left(\boldsymbol M(A_{i_ti_{t+1}})-\boldsymbol\mu_{ij}\right)]^{\lambda_t}\right]^2\\ \nonumber
&=&\sum_{\substack{i_1<i_2<\dots<i_k\\ j_1<j_2<\dots<j_k\\ (i_t,i_{t+1})=(j_t,j_{t+1})\ if\ \lambda_t=1}}\mathbb{E}\left[\prod_{t=1}^k[\boldsymbol 1^T\left(\boldsymbol M(A_{i_ti_{t+1}})-\boldsymbol\mu_{ij}\right)]^{2\lambda_t}\right]=O\left(n^{2k-s}\right),
\end{eqnarray}
where $s\geq |\boldsymbol\lambda|+1$. Then
\begin{eqnarray}\nonumber
&&\sum_{i_1<i_2<\dots<i_k}\sum_{|\boldsymbol\lambda|>0}\prod_{t=1}^k[\boldsymbol 1^T\left(\boldsymbol M(A_{i_ti_{t+1}})-\boldsymbol\mu\right)]^{\lambda_t}[\boldsymbol 1^T(\boldsymbol\mu-\overline{\boldsymbol M}( \boldsymbol A))]^{1-\lambda_t}\\ \nonumber
&=&(1+o_p(1))\prod_{t=1}^k[\boldsymbol 1^T(\boldsymbol\mu-\overline{\boldsymbol M}( \boldsymbol A))]^{1-\lambda_t}\sum_{i_1<i_2<\dots<i_k}\sum_{|\boldsymbol\lambda|<k}\prod_{t=1}^k[\boldsymbol 1^T\left(\boldsymbol M(A_{i_ti_{t+1}})-\boldsymbol\mu_{ij}\right)]^{\lambda_t}\\ \nonumber
&=&(1+o_p(1))O_p\left(\frac{1}{n^{k-|\boldsymbol\lambda|}}\right)O_p\left(2^k\sqrt{n^{2k-s}}\right)\\ \nonumber
&=&(1+o_p(1))O_p\left(2^kC^kn^{|\boldsymbol\lambda|-\frac{s}{2}}\right)\\ \nonumber
&=&(1+o_p(1))O_p\left(2^kC^kn^{\frac{k-1}{2}}\right)\\ \label{power5}
&=&o_p\left(n^{\frac{k}{2}}\left(\boldsymbol 1^T\boldsymbol\mu_{\boldsymbol d}\right)^k\right).
\end{eqnarray}
Then the proof is complete by (\ref{power1}),  (\ref{power3}), (\ref{power2}), (\ref{power4}) and (\ref{power5}).


\begin{thebibliography}{9}


\bibitem{A14}
Aicher, C. (2014).
The Weighted Stochastic Block Model.
\textit{Applied Mathematics Graduate Theses \& Dissertations}, 50.

\bibitem{AJC15}
Aicher, C., Jacob, A. and Clauset, A.(2015).
Learning Latent Block Structure in Weighted Networks.
\textit{Journal of Complex Networks}, 3, 221-248. 

\bibitem{A17} Abbe, E. (2017).
Community detection and stochastic block models: recent developments. 
\textit{Journal of Machine Learning Research}, \textbf{18}, 1-86.


\bibitem{AV14}
Arias-Castro, E. and N. Verzelen. 2014.
Community detection in dense random networks.
\textit{Annals of Statistics}, 42, 3: 940-969.




\bibitem{ALS18}
Ahn, K., Lee, K. and Suh, C.(2018).
Hypergraph spectral clustering in the weighted stochastic block model.
\textit{IEEE Journal of Selected Topics in Signal Processing}. 12, 959-974.




\bibitem{AS17}Abbe, E. and Sandon, C. (2017).
Proof of the achievability conjectures
for the general stochastic block model.
\textit{Communications on Pure and Applied Mathematics}, \textbf{71(7)}, 1334-1406.

\bibitem{ABH16}Abbe, E., Banderira, A. and Hall, G.(2016).
Exact Recovery in the Stochastic Block Model.
\textit{IEEE transactions on information theory}, 62(1) 471-487.



\bibitem{ACB13} Amini, A., Chen, A. and Bickel, P. (2013).
Pseudo-likelihood methods for community detection in large sparse networks.
\textit{Annals of Statistics}, \textbf{41(4)}, 2097-2122.





\bibitem{B18}Banerjee, D. (2018). 
Contiguity and non-reconstruction results for planted partition models: the dense case.
\textit{Electronic Journal of Probability}, \textbf{23}, 1-28.





\bibitem{BM17} Banerjee, D. and Ma, Z. (2017).
Optimal hypothesis testing for stochastic block models with growing degrees. 
\url{https://arxiv.org/pdf/1705.05305.pdf}.  




\bibitem{BS16} Bickel, P. J. and Sarkar, P. (2016). Hypothesis testing for automated community detection in networks.
\textit{Journal of Royal Statistical Society, Series B}, \textbf{78}, 253-273.


\bibitem{BDER16}
Bubeck,S., Ding, J., Eldan,R. and Rácz, M.(2016)
Testing for high‐dimensional geometry in random graphs. \textit{Random Structures \& Algorithms}, 49(3),503-532.

\bibitem{CY06}Chen, J. and Yuan,B. (2006).
Detecting functional modules in the yeast proteinprotein interaction
network. \textit{Bioinformatics}, \textbf{22(18)}, 2283-2290.


\bibitem{CPV07}
Colizza, V., Pastor-Satorras, R. and Vespignani, A.(2007).
Reaction–diffusion processes and metapopulation models in heterogeneous networks.
\textit{Nature Phys}, 3, 276–282.

\bibitem{CS05}
 Canu, S. and Smola, A. J.(2005). Kernel methods and the exponential family. Neurocomputing, 69:
714–720.










\bibitem{F10} Fortunato,S. (2010). Community detection in graphs. \textit{Physics Reports}, \textbf{486 (3-5)},
75-174.







\bibitem{GL17a} Gao, C. and Lafferty, J. (2017a).
Testing for global network structure using small subgraph statistics.
\url{https://arxiv.org/pdf/1710.00862.pdf}







\bibitem{JKL18}
Jin,J., Ke, Z. and Luo S.(2018). Network global testing by counting graphlets. \textit{International conference on machine learning}, 2333-2341.


\bibitem{L16} Lei, J. (2016).
A goodness-of-fit test for stochastic block models.
\textit{Annals of Statistics}, \textbf{44}, 401-424. 




\bibitem{LR15}Lei, J. and Rinaldo, A. (2015).
Consistency of spectral clustering in stochastic block models. \textit{The Annals of Statistics}, \textbf{43(1)}, 215-237.


\bibitem{LWC15}
Lu, Z., Wen, Y. and Cao, G.(2015).
Community detection in weighted networks: algorithms and applications.
\textit{IEEE Transactions on Parallel and Distributed Systems}, 26(11): 2916-2926.




\bibitem{MNS15} Mossel, E., Neeman, J. and Sly, A. (2015).
Reconstruction and estimation in the planted partition model.
\textit{Probability Theory and Related Fields}, \textbf{162}, 431-461.


\bibitem{MS16} Montanari, A. and Sen, S. (2016). Semidefinite programs on sparse random graphs
and their application to community detection. 
\textit{STOC '16 Proceedings of the forty-eighth annual ACM symposium on Theory of Computing},
 814-827. 

\bibitem{N01}Newman, M. (2001).
Scientific collaboration networks. I. Network construction and fundamental results.
\textit{Physical Review E}, \textbf{64}, 016-131.

\bibitem{NBB17}
C. Nicolini, C. Bordier, and A. Bifone. Community detection in weighted brain connectivity networks beyond
the resolution limit. \textit{Phys. Rev. E}, 146:28–39, 2017.


\bibitem{RA15}
R. Rossi and N. Ahmed. (2015).The Network Data Repository with Interactive Graph Analytics and Visualization, \textit{https://networkrepository.com}.


\bibitem{RS10}
M. Rubinov and O. Sporns. Complex network measures of brain connectivity: Uses and interpretations.
\textit{NeuroImage}, 52:1059–1069, 2010.


\bibitem{SFGHK17}
Sriperumbudur, B., Fukumizu, K., Gretton, A., Hyvarinen, A. adn Kumar,R.(2017). Density estimation in infinite dimensional exponential families,\textit{Journal of Machine Learning Research}, 18,1-59.


\bibitem{TB11}
Thomas, A. C. and Blitzstein, J. K. (2011). Valued ties tell fewer lies: Why not to dichotomize network
edges with thresholds. 
\url{arXiv:1101.0788}

\bibitem{T18}
Tokuda, T. (2018). Statistical test for
detecting community structure in real-valued edgeweighted graphs. \textit{PLoS ONE} 13(3): e0194079.

\bibitem{VA15}
Verzelen, N., and E. Arias-Castro. 2015.
Community detection in sparse random networks. 
\textit{Ann. Appl. Probab.} 25,6:3465--3510. 



\bibitem{XJL20}
Xu,M., Jog, V. and Loh,P.(2020).
Optimal rates for community estimation in the weighted stochastic block model.
\textit{Annals of Statistics}, 48(1):183-204.

\bibitem{YFS21}
Yuan, M., Feng,Y. and Shang, Z.(2022).
A likelihood-ratio type test for stochastic block models with bounded degrees. \textit{Journal of Statistical Planning and Inference}, 219:98-119.



\bibitem{YYS22} Yuan, M., Yang, F. and Shang, Z. (2022). Hypothesis testing in sparse weighted stochastic block model. \textit{Statistical Papers}, accepted.





\bibitem{YLFS22} Yuan, M., Liu, R., Feng, Y. and Shang, Z.
(2022). Testing community structures for hypergraphs. 
\textit{Annals of Statistics}, 50(1): 147-169.

\bibitem{YN20}
Yuan, M. and Nan, Y. (2020). Test dense subgraphs in sparse uniform
hypergraph. \textit{Communications in Statistics - Theory and Methods}, to appear.

\bibitem{YS21} Yuan, M. and Shang, Z. (2021). Sharp Detection Boundaries on Testing Dense Subhypergraph.\textit{Bernoulli}, to appear.



\bibitem{ZLZ11} Zhao, Y.,  Levina, E. and  Zhu., J.(2011).
Community extraction for social networks.
\textit{Proc. Natn. Acad. Sci. USA}, \textbf{108}, 7321-7326.

\bibitem{ZLZ12} Zhao, Y.,   Levina, E. and  Zhu, J. (2012).
Consistency of community detection in networks under degree-corrected stochastic block models.
\textit{Annals of Statistics}, \textbf{40}, 2266-2292.























\end{thebibliography}
\end{document}